\documentclass[11pt]{article}
\usepackage[top=1in, bottom=1in, left=1in, right=1in]{geometry}
\usepackage{tikz}
\usepackage{amsmath,amsfonts,amssymb,amscd,amsthm,mathrsfs,float}
\usepackage{extarrows,textcomp}
\usepackage{graphicx,subfigure}
\usepackage{diagbox,multirow}
\usepackage{algorithm,algorithmic}
\usepackage{caption}
\usepackage{enumerate}
\usepackage{bm}
\usepackage{lineno}
\usepackage{xcolor}
\usepackage{makeidx,hyperref}
\usepackage{cite}

\newtheorem{remark}{Remark}[section]
\numberwithin{equation}{section}

\DeclareMathOperator*{\argmin}{arg\,min}
\DeclareMathOperator*{\Div}{div}

\newcommand{\dd}{\,{\rm d}}

\makeatletter
\newcommand*{\extendadd}{
  \mathbin{
    \mathpalette\extend@add{}
  }
}
\newcommand*{\extend@add}[2]{
  \ooalign{
    $\m@th#1\leftrightarrow$%
    \vphantom{$\m@th#1\updownarrow$}
    \cr
    \hfil$\m@th#1\updownarrow$\hfil
  }
}
\makeatother

\begin{document}

\title{An Augmented Lagrangian Deep Learning Method for Variational Problems with Essential Boundary Conditions}

\author{
Jianguo Huang ({\tt jghuang@sjtu.edu.cn})
\vspace{0.1in}\\
School of Mathematical Sciences, and MOE-LSC,
Shanghai Jiao Tong University, Shanghai, China
\vspace{0.1in}\\
Haoqin Wang ({\tt wanghaoqin@sjtu.edu.cn})
\vspace{0.1in}\\
School of Mathematical Sciences, and MOE-LSC,
Shanghai Jiao Tong University, Shanghai, China
\vspace{0.1in}\\
Tao Zhou \footnote{Corresponding author.} ({\tt tzhou@lsec.cc.ac.cn})
\vspace{0.1in}\\
LSEC, Institute of Computational Mathematics and Scientific/Engineering Computing, \\
Academy of Mathematics and Systems Science, Chinese Academy of Sciences, Beijing, China
}

\date{}
\maketitle

\begin{abstract}
This paper is concerned with a novel deep learning method for variational problems with essential boundary conditions. To this end, we first reformulate the original problem into a minimax problem corresponding to a feasible augmented Lagrangian, which can be solved by the augmented Lagrangian method in an infinite dimensional setting. Based on this, by expressing the primal and dual variables with two individual deep neural network functions, we present an augmented Lagrangian deep learning method for which the parameters are trained by the stochastic optimization method together with a projection technique. Compared to the traditional penalty method, the new method admits two main advantages: i) the choice of the penalty parameter is flexible and robust, and ii) the numerical solution is more accurate in the same magnitude of computational cost. As typical applications, we apply the new approach to solve elliptic problems and (nonlinear) eigenvalue problems with essential boundary conditions, and numerical experiments are presented to show the effectiveness of the new method.
\end{abstract}

\\


\section{Introduction}\label{sec: intro}

Variational problems play important roles in various industrial and engineering applications, with typical examples including partial differential equations (PDEs) and eigenvalue problems. Many classical numerical  methods have been developed for such problems, e.g., the finite difference method, the spectral method, and the finite element method. The first two methods are generally used for solving problems over regular domains while the latter one is particularly suitable for problems in irregular domains \cite{Ciarlet1978,BrennerScott2008}. In recent years, deep learning based techniques have been widely used to solve a variety of variational problems \cite{RaissiPerdikarisKarniadakis2019,ZhangGuoKarniadakis2020,PINNRev,EYu2018,SirignanoSpiliopoulos2018,ZangBaoYeZhou2020,HuangWangYang2020,ChenHuangWangYang2020,LiuCaiXu2020,EHanJentzen2017,HanJentzenE2018}. Historically, related studies can date back to the 1990s~\cite{LeeKang1990,DissanayakePhan-Thien1994}.  We also refer the reader to \cite{ERev} and the references therein for a comprehensive review on machine learning from the perspective of computational mathematics. For such kind of methods, deep neural networks (DNNs) are exploited to parameterize the PDE solutions and appropriate parameters are identified by minimizing an optimization problem formulated from the PDEs. The most significant feature of those methods is that they are mesh-free, and their approximation capacity has been well studied in recent years \cite{HornikStinchcombeWhite1989,Barron1993,EMaWu2019,HeLiXuZheng2020,ShenYangZhang2020}.

For variational problems with natural boundary conditions, one doesn't need to impose these conditions on the admissible functions \cite{EYu2018}, so that the DNNs can easily be used for approximation. However, for variational problems with essential boundary conditions, these conditions should be imposed on the admissible functions, and this gives rise to a significant difficulty since one cannot enforce the boundary condition in a simple way even at the interpolation nodes for a neural network function. It is worth noting that even in the context of finite element methods, this is also a very tough issue. In fact, one has to use Nitsche's trick \cite{Nitsche}, developed further by Stenberg \cite{Stenberg1995}, to handle this issue. As far as we know, there are two main strategies to overcome the bottleneck in deep learning framework:
\begin{itemize}
\item
The first strategy is to construct neural network functions that satisfy the essential boundary conditions exactly. For instance, if the boundary condition is given by $u=g$ on the boundary $\Gamma$, then we construct the approximate function by
\begin{equation}
\label{b-condition}
	\phi(\bm{x};\bm{\theta})=\ell(\bm{x})\psi(\bm{x};\bm{\theta})+\bar{g}(\bm{x}),
\end{equation}
where $\ell(\bm{x})$ is a known function such that on $\Gamma$ it holds $\ell(\bm{x})=0$,  $\bar{g}$ is the extension of $g$ to the whole domain, and $\psi(\bm{x};\bm{\theta})$ is another neural network function that is used to approximate the solution in the domain. The main limitation of this approach is that for problems with complex (non-regular) {domains}, it is in general not easy to find explicit functions $\ell$ and $\bar{g}.$ For details, one can refer to \cite{BergNystrom2018} and references therein. It is worth mentioning that based on the formulation \eqref{b-condition}, one can also introduce an additional neural network function on the boundary $\Gamma$ to approximate $g$ by a least squares approach \cite{ShengYang2021}.

\item Another strategy is the penalty method, where a penalty term (with a penalty parameter $\beta$) is included into the objective functional to enforce the boundary condition \cite{RaissiPerdikarisKarniadakis2019,EYu2018,SirignanoSpiliopoulos2018,ZangBaoYeZhou2020}. This method is easy to implement. Theoretically, the penalty parameter $\beta$ should be chosen large enough, however, this may make the optimization problem become ill-conditioned \cite{NocedalWright2006}. We also mention the deep Nitsche method proposed in \cite{LiaoMing2021}, where Nitsche's variational formula is used for the second order elliptic problems to avoid a large penalty parameter.

\end{itemize}

In this work, we intend to present an augmented Lagrangian deep learning (ALDL) method to handle variational problems with essential boundary conditions. For this purpose, we shall first rewrite the original problem as a minimax problem associated with a feasible augmented Lagrangian, which can be solved by the augmented Lagrangian method in an infinite dimensional setting \cite{GlowinskiLeTallec1989,NocedalWright2006}. We then express the primal and dual variables with two individual DNN functions, respectively, and train the associated network parameters with the stochastic optimization method (based on the augmented Lagrangian method). It is worth noting that we require to solve a least square problem in order to update the dual variable, and this step can be viewed as a nonlinear projection in the corresponding parameter space. As typical applications, we apply the ALDL method to solve elliptic problems and eigenvalue problems with essential boundary conditions.  Numerical results indicate that the ALDL method admits two main advantages compared to the penalty method:
\begin{itemize}
\item The choice of the penalty parameter is flexible and robust.
\item The numerical solution is more accurate in the same magnitude of computational cost.
\end{itemize}
The rest of this paper is organized as follows.
In Section \ref{sec: prob}, we introduce the variational problem and its minimax formulation. As typical cases, an elliptic PDE and eigenvalue problems are presented. In Section \ref{sec: alg}, we recall the augmented Lagrangian method in an infinite dimensional setting, and then propose the augmented Lagrangian deep learning method. In Section \ref{sec: numer}, numerical examples are reported to show the performance of the proposed method. Finally, we provide with some concluding remarks in Section 5.

\section{The variational problem and its primal-dual formulation}\label{sec: prob}

To begin, we first introduce some notations for later uses. For a real Hilbert space $V$ equipped with a norm $\| \cdot \|_V$, we denote by $\langle \cdot,\cdot \rangle_{V}$  the induced inner product over $V$. We use the standard symbols and notations for Sobolev space and their norms or semi-norms and refer the reader to the reference \cite{Adams1975} for details. Let $\Omega\subset \mathbb{R}^d$ be a bounded domain with Lipschitz boundary, where $d$ is a natural number. We also denote by $\Gamma$ its boundary and $\bar{\Omega}$ the closure of $\Omega$, respectively. We let $B$ be a bounded linear operator from $V$ to another Hilbert space $W$ and a typical example in our variational problems is $W=L^2(\Gamma)$.

Throughout this paper, we consider the following variational problem
\begin{equation}
\label{VP}
\min_{v \in V_g} J (v),
\end{equation}
where $J(v)$ is a nonlinear functional over $V_g$ and
\[
V_g =\{v\in V: {B} v=g \ \mbox{on} \ \Gamma\}.
\]

\subsection{Primal-dual formulation}\label{subsec:pd}

To deal with the constraint in the admissible set $V_g$, the augmented Lagrangian method \cite{GlowinskiLeTallec1989} suggests to consider an augmented Lagrangian function as following
\begin{equation}
\label{AL}
\mathcal{L}_{\beta}(v,\mu) =
J(v) - \langle \mu, {B} v - g \rangle_{W} + \frac{\beta}{2}\|{B} v - g\|_W^2,
\end{equation}
where $\mu\in W$ is a Lagrange multiplier function (the dual variable) and $\beta$ is a positive constant. Then we obtain the following minimax problem:
\begin{equation}
\label{minmax}
 \min_{v\in V}\max_{\mu \in W}\mathcal{L}_{\beta}(v,\mu).
\end{equation}
By a direct manipulation, we have
\[
  \max_{\mu \in W}\mathcal{L}_{\beta}(v,\mu) =
  \begin{cases}
    J(v),\ &v\in V_g,
    \\
    +\infty, \ &v\notin V_g.
  \end{cases}
\]
Consequently, the variational problem \eqref{VP} is equivalent to the minimax problem \eqref{minmax}. In other words, if $(u,\lambda)\in V\times W$ is a solution of the problem \eqref{minmax}, then $u$ is a solution of the problem \eqref{VP}. On the contrary, if $u$ is a solution of the problem \eqref{VP}, there exists a function $\lambda\in W$ such that $(u,\lambda)$ is a solution of the problem \eqref{minmax}. Furthermore, we assume that $(u,\lambda)$ is a saddle point of the Lagrangian $\mathcal{L}_{\beta}(\cdot,\cdot)$, i.e.,
\begin{equation}
\label{saddle-point}
\min_{v \in V} \max_{\mu \in W} \mathcal{L}_{\beta}(v,\mu)=\mathcal{L}_{\beta}(u,\lambda)=\max_{\mu \in W}\min_{v \in V}\mathcal{L}_{\beta}(v,\mu).
\end{equation}
Notice that it is rather difficult to show the existence of a saddle point of a general functional and one important technique is the Ky Fan-Sion theorem \cite{Cea1978}. If we write
\[
\mathcal{F}_{\beta}(\mu)=\min_{v \in V}\mathcal{L}_{\beta}(v,\mu),
\]
then one can turn to solve the following dual problem
\begin{equation}
\label{dual-problem}
\max_{\mu \in W} \mathcal{F}_{\beta}(\mu).
\end{equation}

\subsection{Some applications}
\label{subsec: apps}

We now present three typical applications of the above primal dual formulation.

\subsubsection{Second order elliptic PDEs}

The first example is the second order elliptic {equation}:
\begin{equation}
\label{elliptic PDE}
\begin{cases}
 -\Div\big(\boldsymbol{A}(\boldsymbol{x})\nabla u(\boldsymbol{x})\big) + c(\boldsymbol{x})u(\boldsymbol{x}) = f(\boldsymbol{x}) & \mbox{\rm in} \ \Omega,
 \\
 u(\boldsymbol{x}) = g(\boldsymbol{x}), &\mbox{\rm on} \ \Gamma,
\end{cases}
\end{equation}
where $\boldsymbol{A}(\boldsymbol{x}) \in C^1(\bar{\Omega})$ is uniformly elliptic and $c(\boldsymbol{x}) \in L^2({\Omega})$ is nonnegative over $\Omega$.

The variational formula of \eqref{elliptic PDE} is
\begin{equation}
\label{Ritz}
\min_{v \in V_g} J(v),
\end{equation}
where
\[
J(v) =
\frac{1}{2}\int_{\Omega} [\boldsymbol{A}(\boldsymbol{x})\nabla v(\boldsymbol{x}) \cdot \nabla v(\boldsymbol{x}) +c (\boldsymbol{x})v^2(\boldsymbol{x}) -2f(\boldsymbol{x})v(\boldsymbol{x})] \dd x,
\]
and
\[
  V_g = \{v\in H^1(\Omega): v = g \ \mbox{a.e. on} \ \Gamma\}.
\]
Based on the abstract setting given in the last subsection, we can rewrite the variational formula \eqref{Ritz} as the following problem
\begin{equation}
\label{PDF-PDE}
\min_{v \in V} \max_{\mu \in W} \mathcal{L}_{\beta}(v,\mu),
\quad
\mathcal{L}_{\beta}(v,\mu)
= J(v) -\int_{\Gamma}\big[\mu(\boldsymbol{x})\big(v(\boldsymbol{x})-g(\boldsymbol{x})\big) \big]\dd x+ \frac{\beta}{2}\int_{\Gamma} \big[v(\boldsymbol{x})-g(\boldsymbol{x})\big]^2 \dd x,
\end{equation}
where $V =H^1(\Omega), W = L^2(\Gamma)$, and $\beta$ is a positive parameter. As shown in \cite{Ciarlet2013}, there exists a unique saddle point for the above problem.

\subsubsection{Linear eigenvalue problems}

The second example is the eigenvalue {problem}. Suppose we want to find the smallest eigenvalue and its eigenfunction for a positive self-adjoint differential operator, e.g.,
\begin{equation}
\label{LinEig}
  \begin{cases}
    - \nabla\cdot\big(p(\boldsymbol{x})\nabla u(\boldsymbol{x})\big)+q(\boldsymbol{x})u(\boldsymbol{x}) = \rho u(\boldsymbol{x}) & \mbox{\rm in} \ \Omega,
    \\
    u(\boldsymbol{x}) =0 \quad  & \mbox{\rm on} \ \Gamma,
  \end{cases}
\end{equation}
where $p(\boldsymbol{x})\in C^1(\bar{\Omega})$ is uniformly elliptic and $q(\boldsymbol{x})\in C(\bar{\Omega})$ is nonnegative over $\Omega$.

The variational formula of {\eqref{LinEig}} (cf. \cite{AmbrosettiMalchiodi2007}) is
\begin{equation}
\label{Rayleigh}
  \min_{v \in V_0} J(v),
\quad
J(v) =
\frac{\int_{\Omega} [p(\boldsymbol{x})\nabla v(\boldsymbol{x}) \cdot \nabla v(\boldsymbol{x}) +q (\boldsymbol{x})v^2(\boldsymbol{x})] \dd x}{\int_{\Omega} v^2(\boldsymbol{x}) \dd x},
\end{equation}
where $V_0=H_0^1(\Omega)$.

Suppose $(\rho,u)$ is the solution of \eqref{LinEig}. Then it is easy to check $(\rho, cu)$ is also a solution of \eqref{LinEig}, where $c$ is a non-zero real number. This motivates us to find the normalized eigenfunction, i.e.,
\begin{equation}
\label{modify Rayleigh}
  \min_{v \in V_0} \tilde{J}(v),
\quad
\tilde{J}(v) = \int_{\Omega} [p(\boldsymbol{x})\nabla \tilde{v}(\boldsymbol{x}) \cdot \nabla \tilde{v}(\boldsymbol{x}) +q (\boldsymbol{x})\tilde{v}^2(\boldsymbol{x})] \dd x,
\end{equation}
where $\tilde{v} = v/\|v\|_{0,\Omega}$.
As mentioned above, we can reformulate \eqref{modify Rayleigh} as the following minimax problem
\begin{equation}
\label{PDF-LinEig}
\min_{v \in V} \max_{\mu \in W} \mathcal{L}_{\beta}(v,\mu),
\quad
\mathcal{L}_{\beta}(v,\mu)
= \tilde{J}(v) - \int_{\Gamma}\mu(\boldsymbol{x})\tilde{v}(\boldsymbol{x}) \dd x + \frac{\beta}{2}\int_{\Gamma} \tilde{v}^2(\boldsymbol{x}) \dd x,
 \end{equation}
where $V=H^1(\Omega)$, $W=L^2(\Gamma)$ and $\beta$ is a positive parameter.

\subsubsection{Nonlinear eigenvalue problems}
The third example is a nonlinear eigenvalue problem:
\begin{equation}
\label{GP}
  \begin{cases}
    -\nabla\cdot\big(\boldsymbol{A}(\boldsymbol{x}) \nabla u(\boldsymbol{x}) \big) + V(\boldsymbol{x})u(\boldsymbol{x}) + u^3(\boldsymbol{x})= \rho u(\boldsymbol{x}) &\quad \mbox{in} \ \Omega,
    \\
    u(\boldsymbol{x}) =0 &\quad \mbox{on} \ \partial\Omega, \\
    \| u \|_{0,\Omega} =1,
  \end{cases}
\end{equation}
where $\boldsymbol{A}(\boldsymbol{x}) \in (L^{\infty}(\Omega))^{d \times d}$ is symmetric and uniformly elliptic, and $V(\boldsymbol{x}) \in L^2(\Omega)$. 

The variational formula of \eqref{GP} (cf. \cite{CancesChakirMaday2010}) is
\begin{equation}\label{GP variational1}
  \min_{\substack{v \in V_0 \\ \|v\|_{0,\Omega}=1}} J(v),
  \quad
  J(v) =
  \frac{1}{2}{\int_{\Omega} \big[ \boldsymbol{A}(\boldsymbol{x}) \nabla v(\boldsymbol{x})\cdot \nabla v(x) + V(\boldsymbol{x})v^2(\boldsymbol{x}) + v^4(\boldsymbol{x}) \big]} \dd x,
\end{equation}
where $V_0 = H_0^1(\Omega)$. According to reference \cite{CancesChakirMaday2010}, the ground state non-negative solution $(\rho, u)$ of \eqref{GP variational1} is unique for $1\leq d \leq3$. Similar to linear eigenvalue problems, we use the normalization technique to relax the constraint $\| v \|_{0,\Omega} = 1$, and then reformulate the variational problem \eqref{GP variational1} as
\begin{equation}
\label{GP variational2}
  \min_{v \in V_0 } \tilde{J}(v),
  \quad
  \tilde{J}(v) =
  \frac{1}{2}{\int_{\Omega} \big[ \boldsymbol{A}(\boldsymbol{x}) \nabla \tilde{v}(\boldsymbol{x})\cdot \nabla \tilde{v}(x) + V(\boldsymbol{x})\tilde{v}^2(\boldsymbol{x}) + \tilde{v}^4(\boldsymbol{x}) \big]} \dd x,
\end{equation}
where $\tilde{v} = v/\|v\|_{0,\Omega}$. As discussed in Section \ref{subsec:pd}, we can rewrite the variational formula \eqref{GP variational2} as the minimax problem
\begin{equation}
\label{PDF-NonEig}
\min_{v \in V} \max_{\mu \in W} \mathcal{L}_{\beta}(v,\mu),
\quad
\mathcal{L}_{\beta}(v,\mu)
= \tilde{J}(v) - \int_{\Gamma}\mu(\boldsymbol{x})\tilde{v}(\boldsymbol{x}) \dd x + \frac{\beta}{2}\int_{\Gamma} \tilde{v}^2(\boldsymbol{x}) \dd x,
\end{equation}
where $V=H^1(\Omega)$, $W=L^2(\Gamma)$ and $\beta$ is a positive parameter.

\section{An augmented Lagrangian deep learning method}
\label{sec: alg}

Before presenting our augmented Lagrangian deep learning method, we  first recall the augmented Lagrangian method for the minimax problem \eqref{minmax} in a Hilbert space setting \cite{GlowinskiLeTallec1989}.

The main idea here is to find the saddle point by solving the dual problem \eqref{dual-problem}. More precisely, one may first fix the primal variable $v_k$ at the $k$th iteration, and update the dual variable to $\mu_{k+1}$. Then we fix the approximate dual variable $\mu_{k+1}$ and find an approximate minimizer $v_{k+1}$. The method is summarized in Algorithm \ref{alg: AL}.

\begin{algorithm}[H]
\small
\caption{The augmented Lagrangian method.}
\label{alg: AL}
\textbf{Input}: The parameter $\beta_0>0$, tolerance $\tau_0 >0$, increase parameter $\alpha \geq 1$, initial guess $v_{0}$, $\mu_0$, max iteration number $Epoch$.
\begin{algorithmic}
\FOR {$k=0,1,\cdots,Epoch$}
  \STATE Update Lagrange multiplier function by
    $\mu_{k+1} = \mu_{k} - \beta_k({B} v_k - g).  \qquad  (\star)$
  \STATE  Fix $\mu_{k+1}$ and find an approximate minimizer $v_{k+1}$ of $\mathcal{L}_{\beta}(v,\mu_{k+1})$ such that $\| \partial_1\mathcal{L}_{\beta}(v_{k+1},\mu_{k+1}) \|\le \tau_k$.
  \STATE Update penalty parameter by $\beta_{k+1}=\alpha \beta_{k}$ .
  \STATE Select new tolerance $\tau_{k+1}$.
\ENDFOR
\end{algorithmic}
\textbf{Output}: $u = v_{Epoch+1}$, $\lambda = \mu_{Epoch+1}$.
\end{algorithm}

We now present some basic concepts of the deep neural networks (DNNs). In this paper, we shall adopt the residual neural network (ResNet) proposed in \cite{HeZhangRenSun2016} to approximate the variational problem. The ResNet can be formulated as follow:
\[
  \boldsymbol{h}_0=\boldsymbol{V}\boldsymbol{x}, \
  \boldsymbol{h}_\ell=\boldsymbol{h}_{\ell-1} + \sigma(\boldsymbol{W}_\ell \boldsymbol{h}_{\ell-1}+\boldsymbol{b}_{\ell}),\ \ell=1,2,\dots,L, \
  {\phi}(\boldsymbol{x};\boldsymbol{\theta})=\boldsymbol{a}^T\boldsymbol{h}_L,
\]
where $\boldsymbol{V}\in \mathbb{R}^{N\times d}$, $\boldsymbol{W}_\ell\in \mathbb{R}^{N\times N}$, $\boldsymbol{b}_\ell\in \mathbb{R}^{N}$ for $\ell=1,\dots,L$, $\boldsymbol{a}\in \mathbb{R}^{N}$. $\sigma(x)$ is a non-linear activation function.
Here, $L$ is the depth of the ResNet, and $N$ is the width of the network.
$\boldsymbol{\theta} =\{ \boldsymbol{V}, \boldsymbol{a}, \boldsymbol{W}_{\ell}, \boldsymbol{b}_{\ell}: 1 \le \ell \le L\}$ denotes the set of all parameters in $\boldsymbol{\phi}$, which uniquely determines the neural network.

To present the augmented Lagrangian deep learning method, we express the primal and dual variables with two individual DNN functions, respectively, i.e.,
\[
u (\boldsymbol{x}) \approx \phi^u(\boldsymbol{x};\boldsymbol{\theta}_{u}), \quad
\lambda (\boldsymbol{x}) \approx \phi^{\lambda}(\boldsymbol{x};\boldsymbol{\theta}_{\lambda}).
\]
To derive the numerical solution, it suffices for us to determine parameters $\boldsymbol{\theta}_{u}$ and $\boldsymbol{\theta}_{\lambda}$. To this end, we may closely  follow the strategies given in Algorithm \ref{alg: AL}. Notice that one can directly follow Algorithm \ref{alg: AL} to update the primal variable, however, for the dual variable, Algorithm \ref{alg: AL} is not directly applicable since we need to update the parameters of the neural networks (not a function itself). In other words, we cannot directly update the parameters $\boldsymbol{\theta}_{\mu}$ of the dual variable $\phi^{\mu}$ through the equation $(\star),$ i.e.,
\[
 \phi^{\mu}_{k+1} \, \Leftarrow \, \phi^{\mu}_k - \beta_k({B} \phi^{v}_k - g).
\]
To overcome this difficulty, we propose to solve the following least squares problem:
\begin{align}
  &\boldsymbol{\theta}^{\mu}_{k+1}
     = \argmin_{\boldsymbol{\theta}_{\nu}} J_{\lambda}(\phi^{\nu}; \phi^{\mu}_k, \phi^{v}_k), \quad \phi^{\nu} := \phi^{\lambda}(\boldsymbol{x};\boldsymbol{\theta}_{\nu}), \ \phi^{\mu}_k := \phi^{\lambda}(\boldsymbol{x};\boldsymbol{\theta}^{\mu}_k), \ \phi^{v}_k := \phi^{u}(\boldsymbol{x};\boldsymbol{\theta}^{v}_k); \nonumber
   \\
  & J_{\lambda}(\phi^{\nu};\phi^{\mu}_k, \phi^{v}_k) = \| \phi^{\nu} - \phi^{\mu}_{k} - \beta_k({B} \phi^{v}_k - g)\|_{L^2(\Gamma)}^2.
\end{align}
Notice that the functionals $\mathcal{L}_{\beta}(v,\mu)$ and $J_{\lambda}(\nu;\mu,v)$ can be expressed as integrals over $\Omega$ or $\Gamma$ in most cases, meaning that sub-optimization problems involved can be solved by means of the stochastic gradient descent method \cite{Bottou2010,KingmaBa2014}.

To sum up the above discussions, we summarize our augmented Lagrangian deep learning method in Algorithm \ref{alg: ALDL}.

\begin{algorithm}[H]
\small
\caption{The augmented Lagrangian deep learning method (ALDL).}
\label{alg: ALDL}
\textbf{Input}: The parameter $\beta_0>0$, the increase parameter $\alpha \geq 1$, the iteration number $Epoch$, the iteration number of model $Epoch_u$, the iteration number of Lagrange multiplier $Epoch_{\lambda}$.
\begin{algorithmic}
\STATE Initialize the network parameters $\boldsymbol{\theta}^v_0$ and $\boldsymbol{\theta}^{\mu}_0$ following the default random initialization of PyTorch. \\
\FOR {$k=0,1,\cdots, Epoch$}
  \STATE Fix the parameter $\boldsymbol{\theta}^{v}_k$ and the parameter $\boldsymbol{\theta}^{\mu}_k$. With the initial guess $\boldsymbol{\theta}^{\mu}_k$, find an approximate solver $\boldsymbol{\theta}^{\mu}_{k+1}$ of $\min_{\boldsymbol{\theta}_{\nu}} J_{\lambda}\big(\phi^\nu;\phi^{\mu}_k,\phi^v_k\big)$ by the Adam method with $Epoch_{\lambda}$ iteration steps.
  \STATE Fix the parameter $\boldsymbol{\theta}^{\mu}_{k+1}$. With the initial guess $\boldsymbol{\theta}^{v}_k$, find an approximate solver $\boldsymbol{\theta}^{v}_{k+1}$ of $\min_{\boldsymbol{\theta}_{u}}\mathcal{L}_{\beta}\big(\phi^u(\cdot \ ;\boldsymbol{\theta}_{u}),\phi^{\mu}_{k+1}\big)$ by the Adam method with $Epoch_{u}$ iteration steps.
  \STATE Update penalty parameter by $\beta_{k+1}=\alpha \beta_{k}$.
\ENDFOR
\end{algorithmic}
\textbf{Output}: $u = \phi^{v}_{Epoch+1}$, $\lambda = \phi^{\mu}_{Epoch+1}$.
\end{algorithm}

\begin{remark}\label{remark}
Another strategy to update the dual variable $\phi^{\mu}$ is solving the minimax problem \eqref{minmax} in the form of DNNs by the stochastic gradient descent ascent method (SGDA) directly. However, how to solve this problem efficiently is still an interesting problem and one can refer {to}\cite{NemirovskiJuditskyLanShapiro08,YoonRyu21,XianHuangZhangHuang2021,HuangWu2021} and the reference therein for more details. We only make some numerical comparisons of the ALDL method and the SGDA method in Section 4.
\end{remark}

Next, we present some details for implementing the ALDL method. For the the elliptic problem \eqref{elliptic PDE}, we have
\begin{align}
  \mathcal{L}_{\beta}(v,\mu)
=
& \frac{1}{2}\int_{\Omega}\big[ \boldsymbol{A}(\boldsymbol{x})\nabla v(\boldsymbol{x}) \cdot \nabla v(\boldsymbol{x}) +c (\boldsymbol{x})v^2(\boldsymbol{x}) -2f(\boldsymbol{x})v(\boldsymbol{x})\big] \dd x \nonumber
\\
&-\int_{\Gamma}\big[\mu(\boldsymbol{x})\big(v(\boldsymbol{x})-g(\boldsymbol{x})\big)\big] \dd x+ \frac{\beta}{2}\int_{\Gamma} \big[v(\boldsymbol{x})-g(\boldsymbol{x})\big]^2 \dd x \nonumber
\\
=
&|\Omega|\mathbb{E}_{\boldsymbol{\xi}}
\left[\frac{1}{2}\boldsymbol{A}(\boldsymbol{\xi})|\nabla v(\boldsymbol{\xi})|^2 + \frac{1}{2}c (\boldsymbol{\xi})v^2(\boldsymbol{\xi})
-f(\boldsymbol{\xi})v(\boldsymbol{\xi})\right] \nonumber
\\
&-|\Gamma|\mathbb{E}_{\boldsymbol{\eta}}
\left[\mu(\boldsymbol{\eta})\big(v(\boldsymbol{\eta})-g(\boldsymbol{\eta})\big)
+ \frac{\beta}{2}\big(v(\boldsymbol{\eta})-g(\boldsymbol{\eta})\big)^2\right],
\label{PDE loss}
\end{align}
where $\boldsymbol{\xi}$ and $\boldsymbol{\eta}$ are random vectors following the uniform distribution over $\Omega$ and $\Gamma$, respectively. $|\Omega|$ and $|\Gamma|$ are the measure of $\Omega$ and $\Gamma$, respectively.

By the definition of {$J_{\lambda}(\nu;\mu,v)$}, we have
\begin{align}
J_{\lambda}(\nu;\mu_k,v_k)
&=
\int_{\Gamma}\left[\nu(\boldsymbol{x})-\mu_k(\boldsymbol{x})-\beta\big(v_k(\boldsymbol{x})-g(\boldsymbol{x})\big)\right]^2 \dd x \nonumber
\\
&=
|\Gamma|\mathbb{E}_{\boldsymbol{\eta}}
\left[\big(\nu(\boldsymbol{\eta})-\mu_k(\boldsymbol{\eta})-\beta(v_k(\boldsymbol{\eta})-g(\boldsymbol{\eta}))\big)^2\right].
\label{PDE lambda}
\end{align}

Upon substituting \eqref{PDE loss} and \eqref{PDE lambda} into Algorithm \ref{alg: ALDL}, one can then obtain the ALDL algorithm.

\vspace{0.1 in}
As for the linear eigenvalue problem, a direct manipulation gives

\begin{align}
  \mathcal{L}_{\beta}(v,\mu)
  &= \frac{\mathbb{E}_{\boldsymbol\xi}\big[p(\boldsymbol{\xi})|\nabla {v}(\boldsymbol{\xi})|^2  +q (\boldsymbol{\xi})v^2(\boldsymbol{\xi})\big]}{\mathbb{E}_{\boldsymbol\zeta}\big[v^2(\boldsymbol{\zeta})\big]}
  - \frac{|\Gamma|\mathbb{E}_{\boldsymbol\eta}\big[\mu(\boldsymbol{\eta}){v}(\boldsymbol{\eta})\big]}{\sqrt{|\Omega|\mathbb{E}_{\boldsymbol\zeta}\big[v^2(\boldsymbol{\zeta})\big]}}
  + \frac{\beta|\Gamma|\mathbb{E}_{\boldsymbol\eta}\big[ v^2(\boldsymbol{\eta})\big]}{2|\Omega|\mathbb{E}_{\boldsymbol\zeta}\big[v^2(\boldsymbol{\zeta})\big] },
  \\
  J_{\lambda}(\nu;\mu_k,\tilde{v}_k)
  &=
  |\Gamma|\mathbb{E}_{\boldsymbol{\eta}}
  \left[\big(\nu(\boldsymbol{\eta})-\mu_k(\boldsymbol{\eta})-\beta \tilde{v}_k(\boldsymbol{\eta})\big)^2\right]
  = |\Gamma|\mathbb{E}_{\boldsymbol{\eta}} \left[\left(\nu(\boldsymbol{\eta})-\mu_k(\boldsymbol{\eta})-\beta \frac{{v}_k(\boldsymbol{\eta})}{\|v_k\|_{0,\Omega}}\right)^2\right] ,
\end{align}
where $\boldsymbol{\xi}$, $\boldsymbol{\zeta}$ are i.i.d. random vectors following the uniform distribution over $\Omega$, and $\boldsymbol{\eta}$ is the random vector following the uniform distribution over $\Gamma$. 
Notice that $\tilde{v}_k$ is fixed in $J_{\lambda}(\nu;\mu_k,\tilde{v}_k)$, so $\|v_k\|_{0,\Omega}$ is a known real number in $J_{\lambda}(\nu;\mu_k,\tilde{v}_k).$ Again, one may obtain the ALDL algorithm for eigenvalue problems by substituting the above formulas into Algorithm \ref{alg: ALDL}.

For the nonlinear eigenvalue problem \eqref{GP}, we have
\begin{align}
  \mathcal{L}_{\beta}(v,\mu)
  &= \frac{\mathbb{E}_{\boldsymbol{\xi}_1}\big[\boldsymbol{A}(\boldsymbol{\xi}_1)|\nabla {v}(\boldsymbol{\xi}_1)|^2
  + V (\boldsymbol{\xi}_1)v^2(\boldsymbol{\xi}_1)\big]}{\mathbb{E}_{\boldsymbol{\xi}_2}\big[v^2(\boldsymbol{\xi}_2)\big]}
  + \frac{\mathbb{E}_{\boldsymbol{\xi}_1}\big[v^4(\boldsymbol{\xi}_1)\big]}{|\Omega|\mathbb{E}_{\boldsymbol{\xi}_2,\boldsymbol{\xi}_3}\big[v^2(\boldsymbol{\xi}_2) v^2(\boldsymbol{\xi}_3)\big]}
   \nonumber  \\
  &- \frac{|\Gamma|\mathbb{E}_{\boldsymbol\eta}\big[\mu(\boldsymbol{\eta}){v}(\boldsymbol{\eta})\big]}{\sqrt{|\Omega|\mathbb{E}_{\boldsymbol{\xi}_2}\big[v^2(\boldsymbol{\xi}_2)\big]}}
  + \frac{\beta|\Gamma|\mathbb{E}_{\boldsymbol\eta}\big[ v^2(\boldsymbol{\eta})\big]}{2|\Omega|\mathbb{E}_{\boldsymbol{\xi}_2}\big[v^2(\boldsymbol{\xi}_2)\big] },
  \\
  J_{\lambda}(\nu;\mu_k,\tilde{v}_k)
  &=
  |\Gamma|\mathbb{E}_{\boldsymbol{\eta}}
  \left[\big(\nu(\boldsymbol{\eta})-\mu_k(\boldsymbol{\eta})-\beta \tilde{v}_k(\boldsymbol{\eta})\big)^2\right]
  = |\Gamma|\mathbb{E}_{\boldsymbol{\eta}} \left[\left(\nu(\boldsymbol{\eta})-\mu_k(\boldsymbol{\eta})-\beta \frac{{v}_k(\boldsymbol{\eta})}{\|v_k\|_{0,\Omega}}\right)^2\right] ,
\end{align}
where $\boldsymbol{\xi}_1$, $\boldsymbol{\xi}_2$, $\boldsymbol{\xi}_3$ are i.i.d. random vectors produced by the uniform distribution over $\Omega$, and $\boldsymbol{\eta}$ is the random vector following the uniform distribution over $\Gamma$.

\section{Numerical experiments}\label{sec: numer}
In this section, we shall present various numerical examples to illustrate the effectiveness of the ALDL method. We shall also perform a numerical comparison between the penalty method (PMDL) and the ALDL method. In all our numerical examples, we shall use the ResNet with width $N=50$ and depth $L=6$ with an activation function $\sigma(x) = \max \{ x, 0 \}^2$ to approximate the solution $u$, which results in 15350 unknown parameters. To reduce the computation cost, we use a ResNet with width $N=50$ and depth $L=2$ with an activation function $\sigma(x) = \max \{ x, 0 \}^2$ as an approximation to the dual variable $\lambda$, which results in 5150 parameters.

We use the Adam optimizer \cite{KingmaBa2014} for training, with a learning rate $\eta = 1\text{e}-3$ for elliptic PDEs and linear eigenvalue problems, and $\eta = 5\text{e}-4$ for nonlinear eigenvalue problems. The learning rate in the optimization is adjusted in an exponentially decaying scheme, where the decaying rate is $0.01^{\frac{1}{50000}}$. The maximum iteration number $Epoch$ is taken as 50000 for the PMDL method, and we set $Epoch=500$, $Epoch_{\lambda}=100$ and $Epoch_u=100$ for the ALDL method. The batch size within the computation domain is 512 for 2d problems and 2048 for 3d problems. While the number of training points on each boundary is set to be 64 {for 2d} problems and is set to be 256 for 3d problems. All numerical experiments are implemented in Python 3.7 using Pytorch 1.3 in an NVIDIA GEFORCE RTX 2080 Ti GPU card. The code can be shared upon request.

For elliptic PDEs, we use $u$ and $u_{dl}$ to denote the exact solution and the deep learning solution, respectively. While for eigenvalue problems, we let $\rho$ and $u$ be the smallest eigenvalue and the corresponding eigenfunction with $\|u\|_0=1$, respectively. Also, we denote by $\rho_{dl}$ and $u_{dl}$ the DNNs solutions.
For a computational domain $\Omega=(0,1)^d$, we divide it into small cubes with a length $h$ uniformly and denote the set of all vertexes as $\Omega_h$, i.e., $\Omega_h = \{ih: 0\leq i \leq N \}^d$ with $N=1/h$. In all experiments, we shall adopt $h=2^{-6}$ to generate the set $\Omega_h$ as the test locations. In order to measure the accuracy of deep learning algorithms, we introduce the discrete maximum norm
\[
  \|v\|_{0,\infty,h}=\max_{\boldsymbol{x}\in E}|v(\boldsymbol{x})|, \quad E\subset\Omega_h.
\]
In addition, we define the absolute error and relative error in the domain $\Omega$ as follows:
\[
  \mathcal{E}_a^{in} = \|u(\boldsymbol{x})-u_{dl}(\boldsymbol{x})\|_{0,\infty,h} ,
  \quad
  \mathcal{E}_r^{in}= \left\| \frac{u(\boldsymbol{x})-u_{dl}(\boldsymbol{x})}{u(\boldsymbol{x})} \right\|_{0,\infty,h},
  \quad
  \boldsymbol{x} \in E = \Omega \cap \Omega_h.
\]
The absolute error and relative error on the boundary $\Gamma$ are defined by
\[
  \mathcal{E}_a^{bd} = \|u(\boldsymbol{x})-u_{dl}(\boldsymbol{x})\|_{0,\infty,h},
  \quad
  \mathcal{E}_r^{bd}= \left\| \frac{u(\boldsymbol{x})-u_{dl}(\boldsymbol{x})}{u(\boldsymbol{x})} \right\|_{0,\infty,h},
  \quad
  \boldsymbol{x} \in E = \Gamma \cap \Omega_h.
\]

\subsection{Elliptic PDEs}

Consider the Poisson equation with the Dirichlet boundary
\[
  \begin{cases}
    - \Delta u = f \quad  &\mbox{\rm in}  \ \Omega,
    \\
    \quad u = g \quad &\mbox{\rm on} \ \Gamma,
  \end{cases}
\]
where $\Omega = (0,1)^d$ for $d=2$ and $3$. We choose appropriate $f$ and $g$ such that the exact solution yields 
\[
u(\bm{x}) = \sum
_{i=1}^d\big(\sin(2\pi x_i)+1.25\big).
\]

\subsubsection{The 2D case}

\begin{figure}[H]
  \begin{minipage}[H]{0.5\linewidth}
    \centering
    \includegraphics[scale = 0.14]{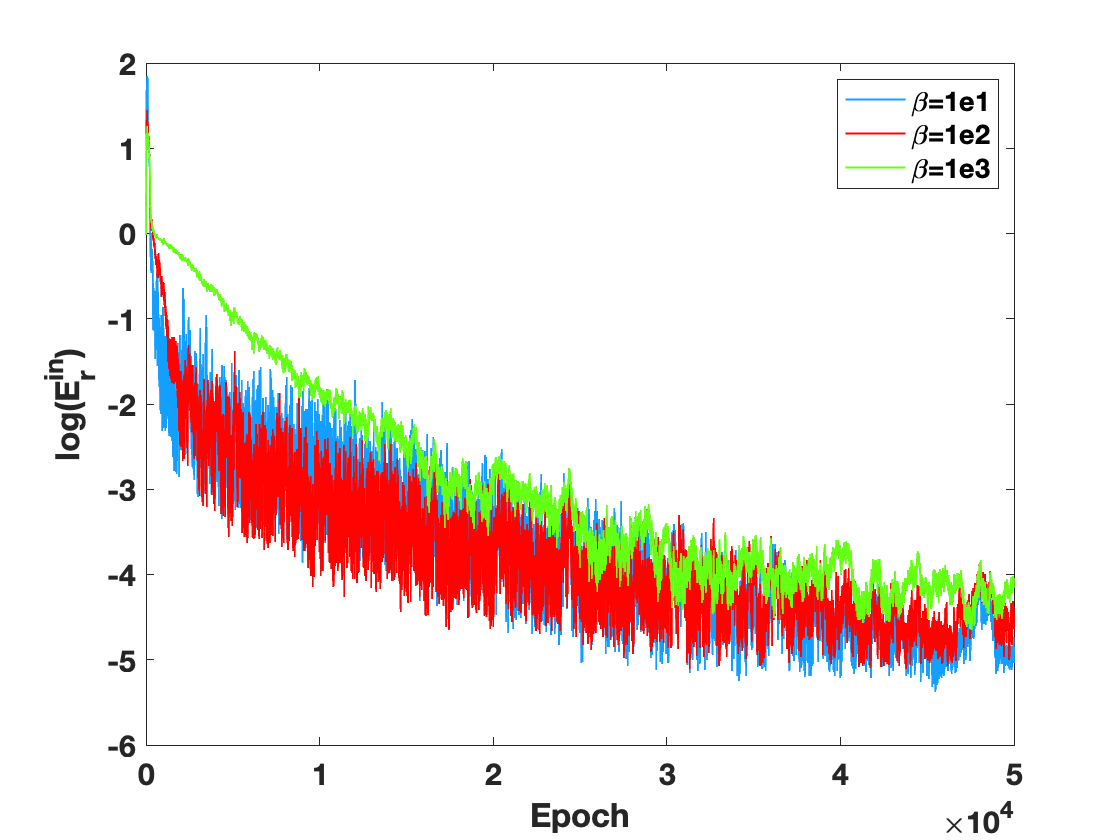}
  \end{minipage}
  \begin{minipage}[H]{0.5\linewidth}
    \centering
    \includegraphics[scale = 0.14]{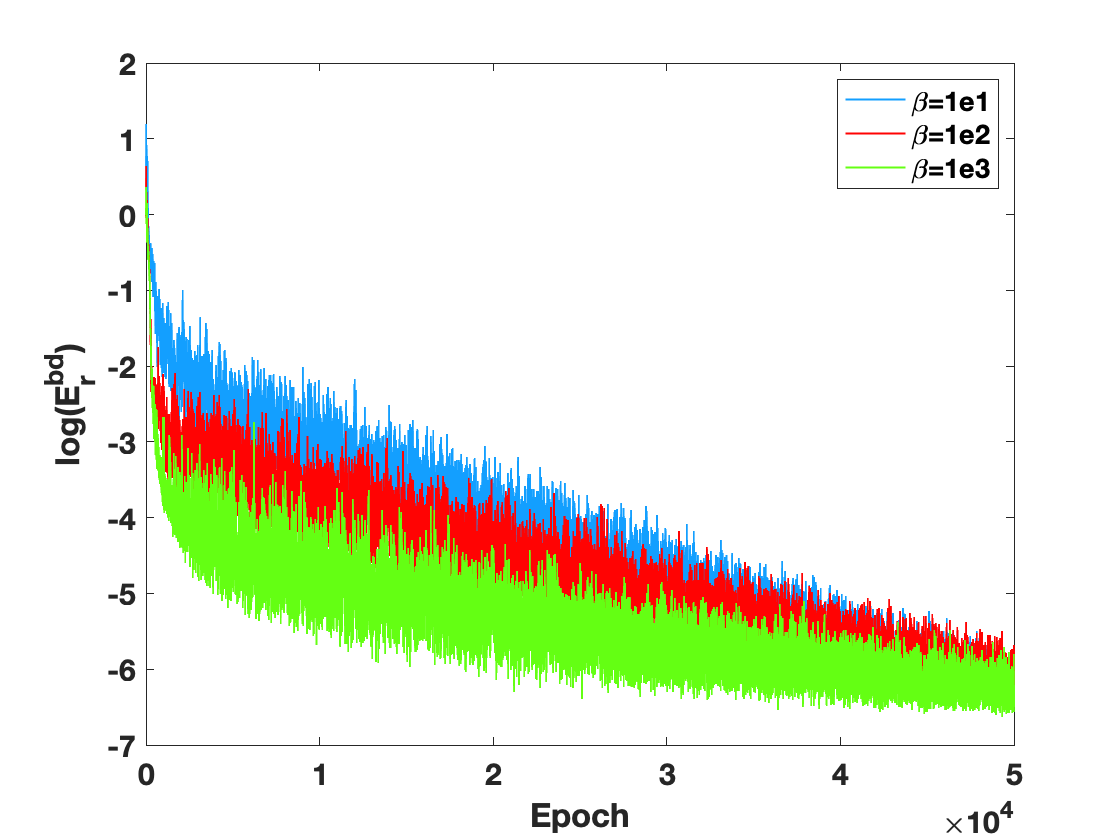}
  \end{minipage}
  \caption{The relative error in the domain (left) and on the boundary (right) with ALDL.}
  \label{fig:2d_pde_al}
\end{figure}

We first solve the above problem by the ALDL method with different parameters $\beta$. The test errors are presented in Figure \ref{fig:2d_pde_al} and in Table \ref{tab:2d_pde_al}.
It can be seen from Figure \ref{fig:2d_pde_al} (left) that the ALDL method is rather robust with respect to $\beta,$ and there is no evident difference on the accuracy (within the domain) and convergent speed when $\beta$ is chosen from {$10$} to $1000$. In other words, one may use a relatively small $\beta$ in practice. By Figure \ref{fig:2d_pde_al} (right), we can also see the ALDL method admits {a} very good accuracy on the boundary.

\begin{table}[H]
\small
\begin{minipage}[H]{0.5\linewidth}
  \centering
  \begin{tabular}{ccc}
    \hline
    $\beta$ &$\mathcal{E}_{r}^{in}$ &$\mathcal{E}_{r}^{bd}$ \\
    \hline
     1e+1  &  7.7843e-3 &  1.9851e-3 \\
     1e+2  &  1.0998e-2 &  1.9644e-3 \\
     1e+3  &  1.6788e-2 &  2.0246e-3 \\
    \hline
  \end{tabular}
  \caption{ The final relative errors of ALDL. }
  \label{tab:2d_pde_al}
\end{minipage}
\begin{minipage}[H]{0.5\linewidth}
  \centering
  \begin{tabular}{ccc}
  \hline
    $\beta$ &$\mathcal{E}_{r}^{in}$ &$\mathcal{E}_{r}^{bd}$ \\
    \hline
     2e+2  &  4.8946e-2 &  2.2827e-2 \\
     2e+3  &  1.2054e-2 &  3.7951e-3 \\
     2e+4  &  4.0631e-1 &  1.5662e-3 \\
    \hline
  \end{tabular}
  \caption{The final relative errors of PMDL. }
  \label{tab:2d_pde_pm}
\end{minipage}
\end{table}

\begin{figure}[H]
  \begin{minipage}[H]{0.5\linewidth}
    \centering
    \includegraphics[scale = 0.14]{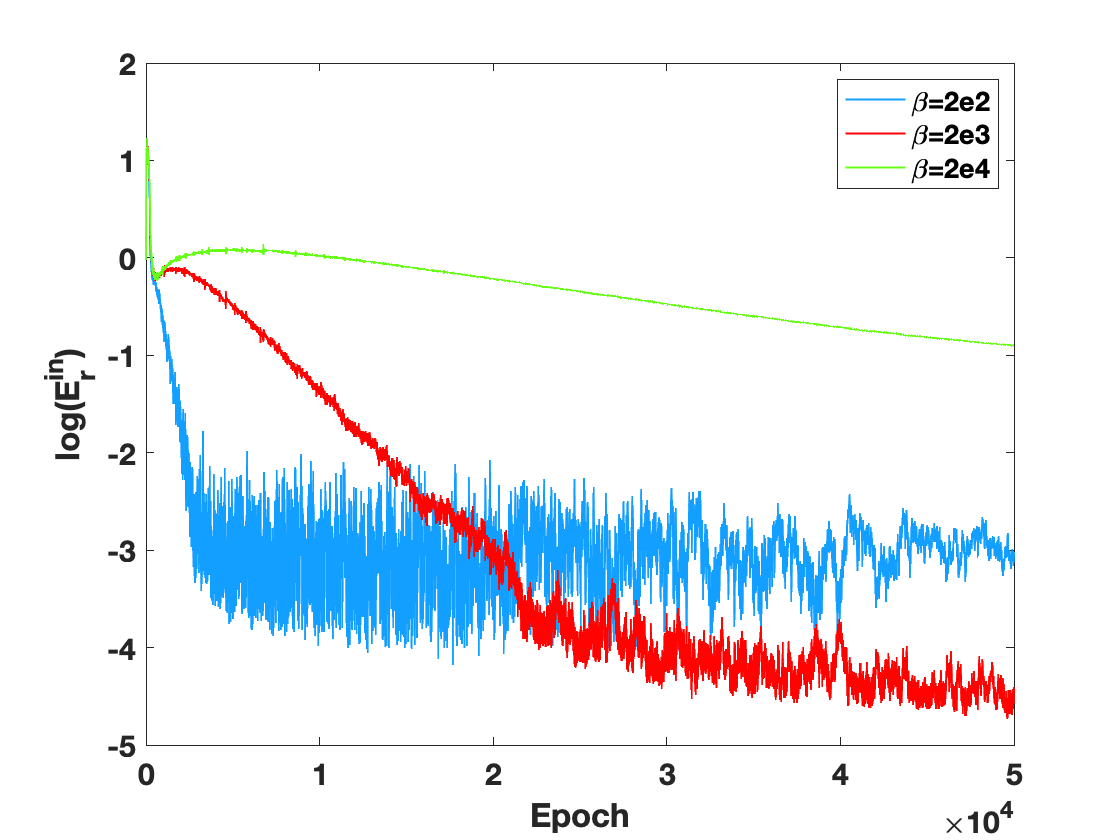}
  \end{minipage}
  \begin{minipage}[H]{0.5\linewidth}
    \centering
    \includegraphics[scale = 0.14]{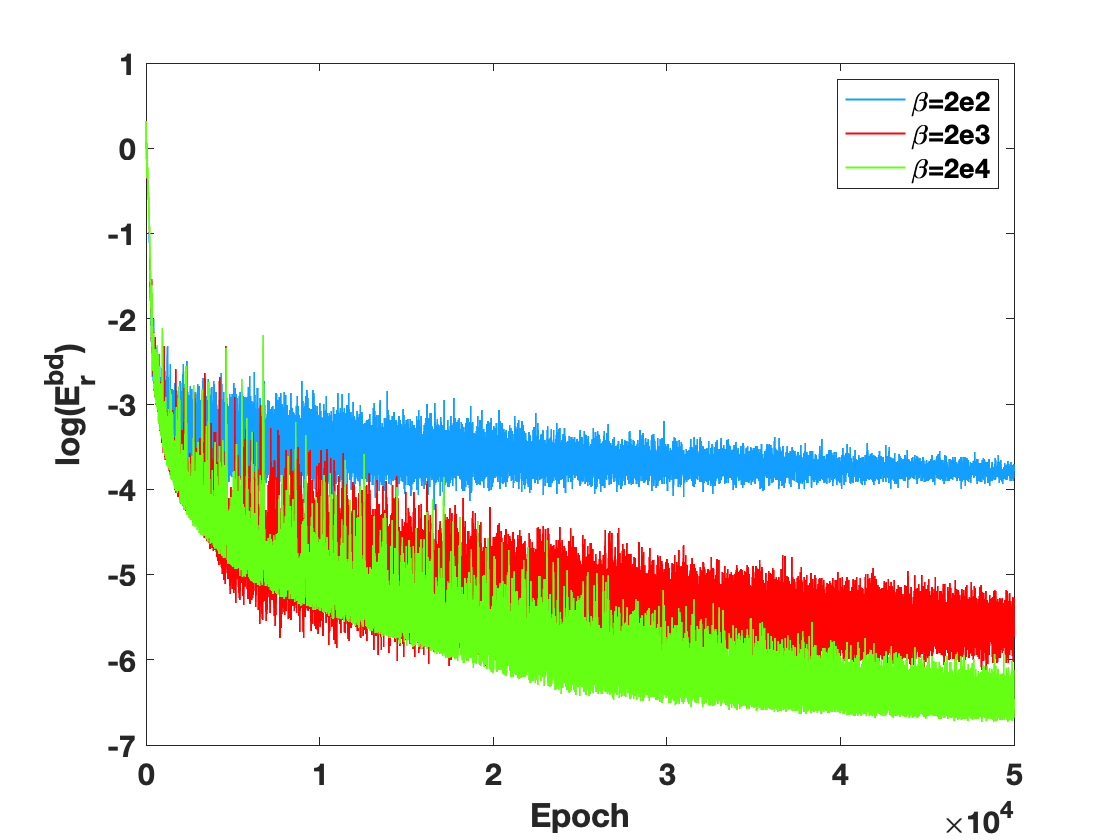}
  \end{minipage}
  \caption{The relative error in the domain (left) and on the boundary (right) of PMDL.}
  \label{fig:2d_pde_pm}
\end{figure}

We also solve the problem with the PMDL method with different parameters $\beta$, and the numerical results are shown in Figure \ref{fig:2d_pde_pm} and Table \ref{tab:2d_pde_pm}. By Figure \ref{fig:2d_pde_pm} (right)  we can observe that one should use a large parameter $\beta$ to get a good accuracy on the boundary. However, by Figure \ref{fig:2d_pde_pm} (left) we can see that the large the parameter $\beta$ is, the slower the converge rate. This may due to the ill-conditioning for a relatively large $\beta$ in the penalty methods \cite{NocedalWright2006}.

As mentioned in Remark \ref{remark}, one can solve the minimax problem \eqref{PDF-PDE} in the form of DNNs by the SGDA method. To compare the numerical performance between the ALDL method and the SGDA method, we list the errors and GPU time for different methods with the same parameters in Table \ref{tab:2d_pde_com}. We learn that the ALDL method and the SGDA method admit a similar accuracy, however, the ALDL method saves about $25\%$ more computational cost than the SGDA method.

\begin{table}[H]
\small
  \centering
  \begin{tabular}{ccccc}
  \hline
    Method &$\beta$ &$\mathcal{E}_{r}^{in}$ &$\mathcal{E}_{r}^{bd}$ &time(s) \\
    \hline
     ALDL &  1e+2  &  1.0998e-2 & 1.9644e-3  & 933.82\\
     SGDA &  1e+2  &  1.0651e-2 & 2.8226e-3  & 1256.22\\
     ALDL &  1e+3  &  1.6788e-2 & 2.0246e-3  & 938.39\\
     SGDA &  1e+3  &  1.7710e-2 & 2.4454e-3  & 1223.07\\
    \hline
  \end{tabular}
  \caption{The comparation of different methods. }
  \label{tab:2d_pde_com}
\end{table}

\subsubsection{The 3D case}

\begin{figure}[H]
  \begin{minipage}[H]{0.5\linewidth}
    \centering
    \includegraphics[scale = 0.14]{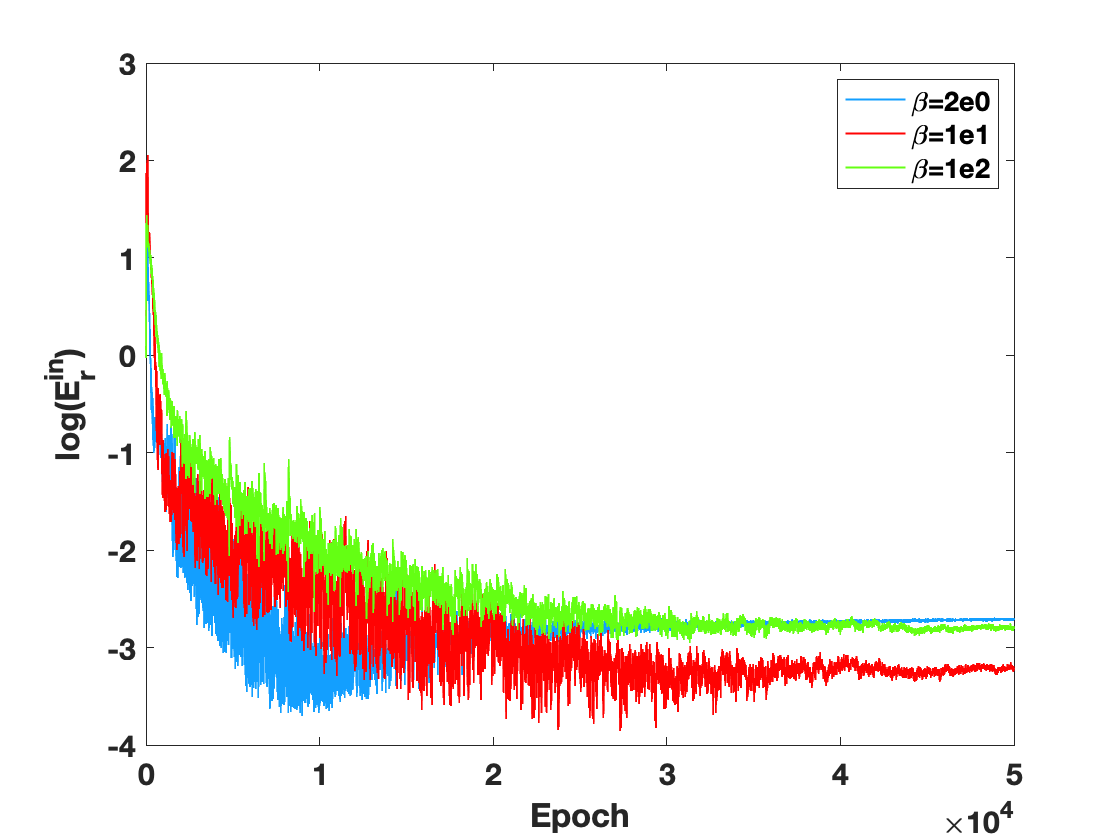}
  \end{minipage}
  \begin{minipage}[H]{0.5\linewidth}
    \centering
    \includegraphics[scale = 0.14]{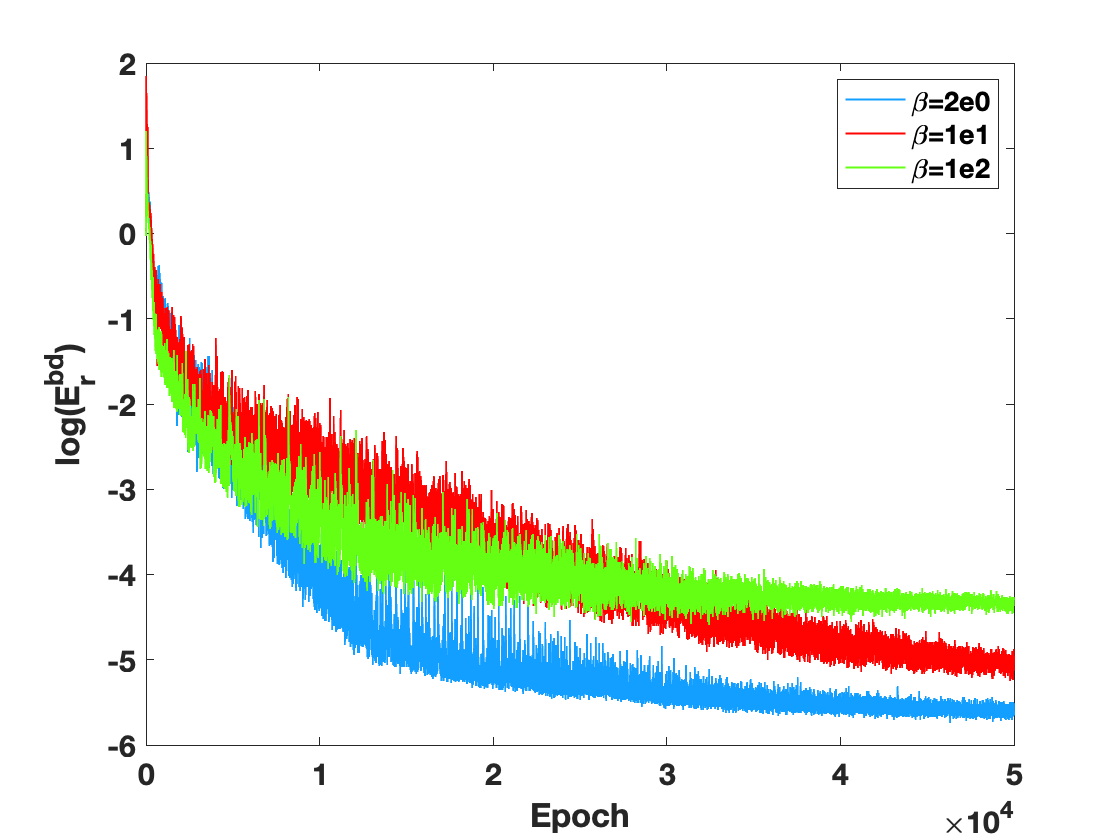}
  \end{minipage}
  \caption{ The relative error in the domain (left) and on the boundary (right) of ALDL.}
  \label{fig:3d_pde_al}
\end{figure}

\begin{figure}[H]
  \begin{minipage}[t]{0.5\linewidth}
    \centering
    \includegraphics[scale = 0.14]{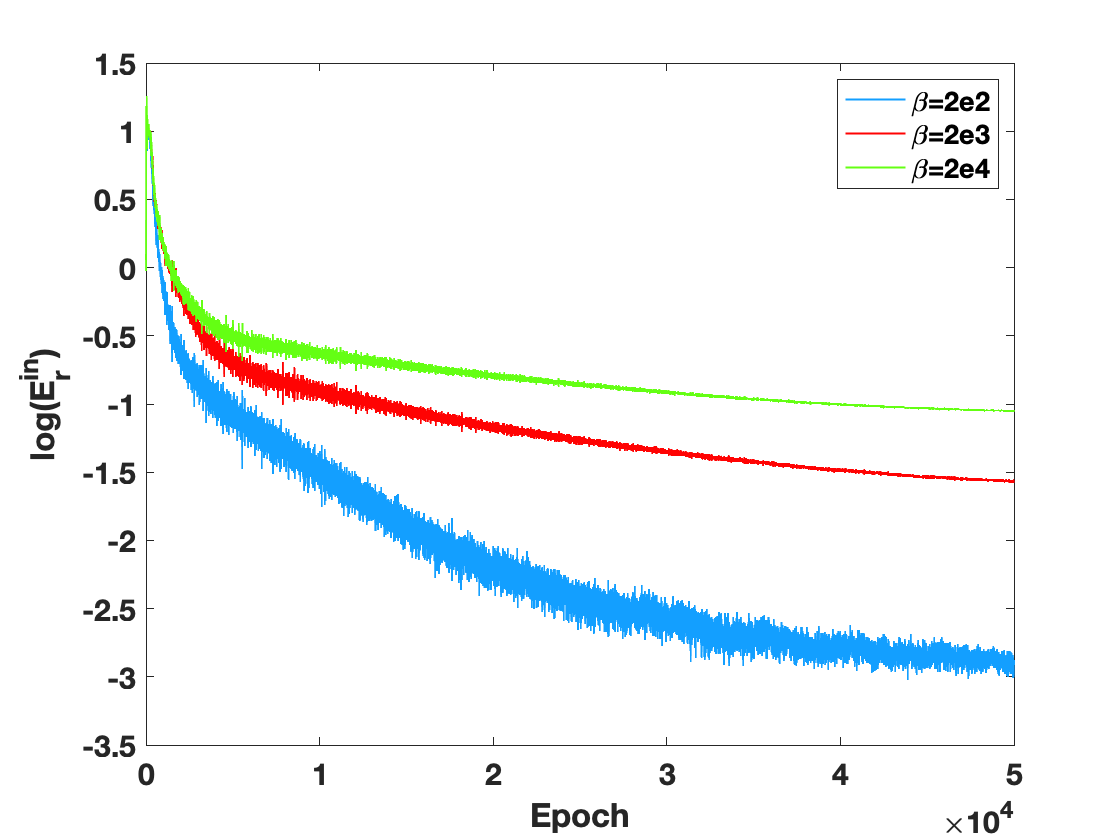}
  \end{minipage}
  \begin{minipage}[t]{0.5\linewidth}
    \centering
    \includegraphics[scale = 0.14]{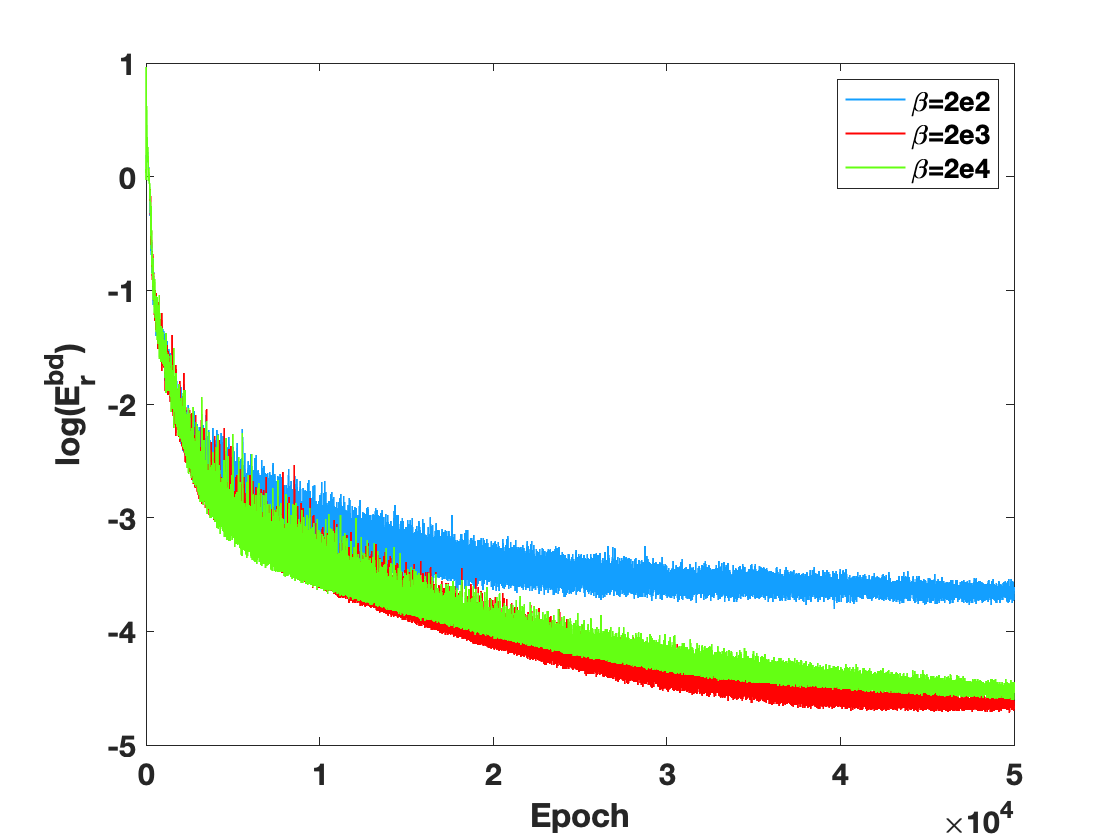}
  \end{minipage}
  \caption{ The relative error in the domain (left) and on the boundary (right) of PMDL.}
  \label{fig:3d_pde_pm}
\end{figure}

We next consider the three dimensional case. Similar plots are shown in  Figure \ref{fig:3d_pde_al} and Figure \ref{fig:3d_pde_pm}, and the final approximation errors are listed in Table \ref{tab:3d_pde_al} and Table \ref{tab:3d_pde_pm}. From those pictures and tables, one can draw similar conclusions as in the two dimensional case. 
{When both methods provide good approximations, the accuracy of the ALDL method is about $2$ times higher than that of the PMDL method. However, in the case of large parameters for PMDL method, the ALDL method outperforms the PMDL method about $50$ times.}

\begin{table}[H]
\small
\begin{minipage}[H]{0.5\linewidth}
  \centering
  \begin{tabular}{ccc}
    \hline
    $\beta$ &$\mathcal{E}_{r}^{in}$ &$\mathcal{E}_{r}^{bd}$ \\
    \hline
     2e+0   &  6.6582e-2 & 3.5099e-3 \\
     1e+1   &  4.0046e-2 & 6.7436e-3 \\
     1e+2   &  5.9976e-2 & 1.3568e-2 \\
    \hline
  \end{tabular}
  \caption{ The final relative errors of ALDL.}
  \label{tab:3d_pde_al}
\end{minipage}
\begin{minipage}[H]{0.5\linewidth}
  \centering
  \begin{tabular}{ccc}
  \hline
    $\beta$ &$\mathcal{E}_{r}^{in}$ &$\mathcal{E}_{r}^{bd}$ \\
    \hline
     2e+2  &  5.6320e-2 &  2.6998e-2 \\
     2e+3  &  2.0876e-1 &  9.2283e-3 \\
     2e+4  &  3.5052e-1 &  1.0640e-2 \\
    \hline
  \end{tabular}
  \caption{ The final relative errors of PMDL. }
  \label{tab:3d_pde_pm}
\end{minipage}
\end{table}

\subsection{Linear eigenvalue problems}
We next consider the following eigenvalue problem
\[
  \begin{cases}
    - \Delta u = \rho u \quad   &\mbox{\rm in}  \ \Omega,
    \\
    \quad u = 0 \quad &\mbox{\rm on} \ \Gamma,
  \end{cases}
\]
where $\Omega = (0,1)^d$ for $d=2$ and $3$. The exact smallest eigenvalue is $d\pi^2$ and the corresponding eigenfunction is 
\[
u(\bm{x}) = \prod_{i=1}^d\sin(\pi(x_i-1)).
\]
Due to the homogeneous Dirichlet boundary, we evaluate different approaches with absolute errors for eigenfunctions (and relative errors for eigenvalues).

\subsubsection{The 2D case}

\begin{figure}[H]
  \begin{minipage}[H]{0.5\linewidth}
    \centering
    \includegraphics[scale = 0.14]{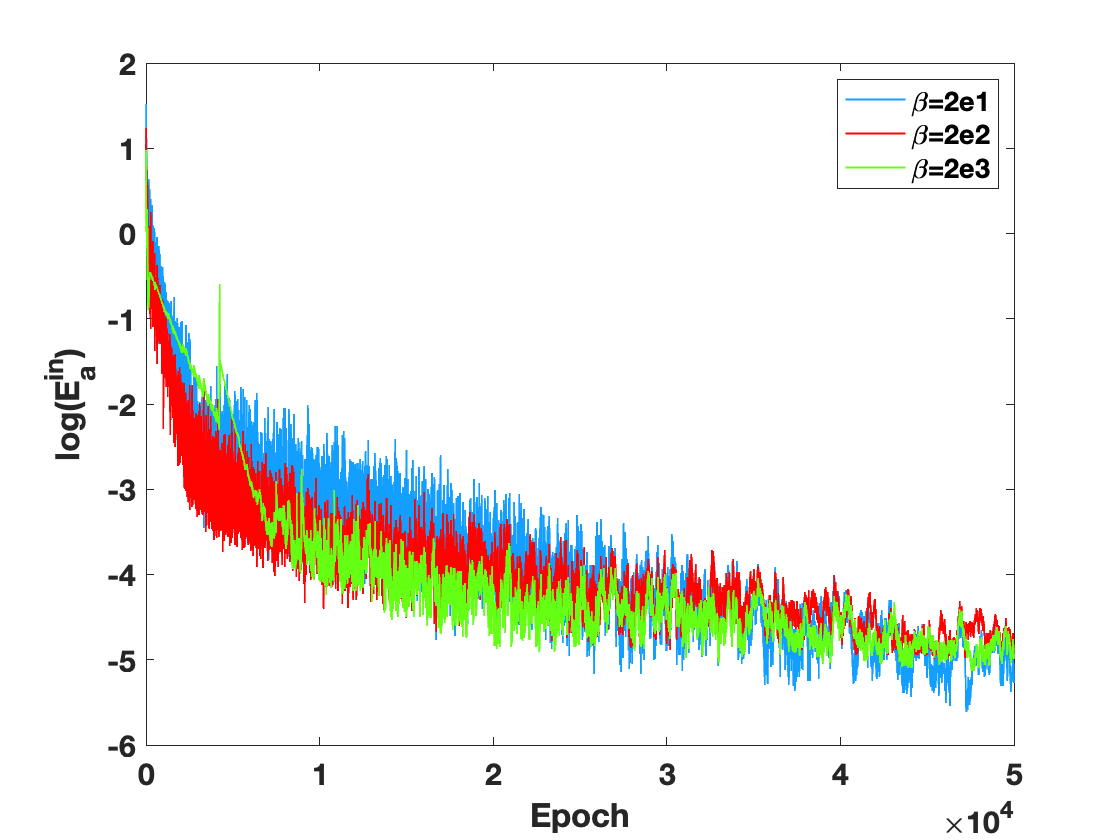}
  \end{minipage}
  \begin{minipage}[H]{0.5\linewidth}
    \centering
    \includegraphics[scale = 0.14]{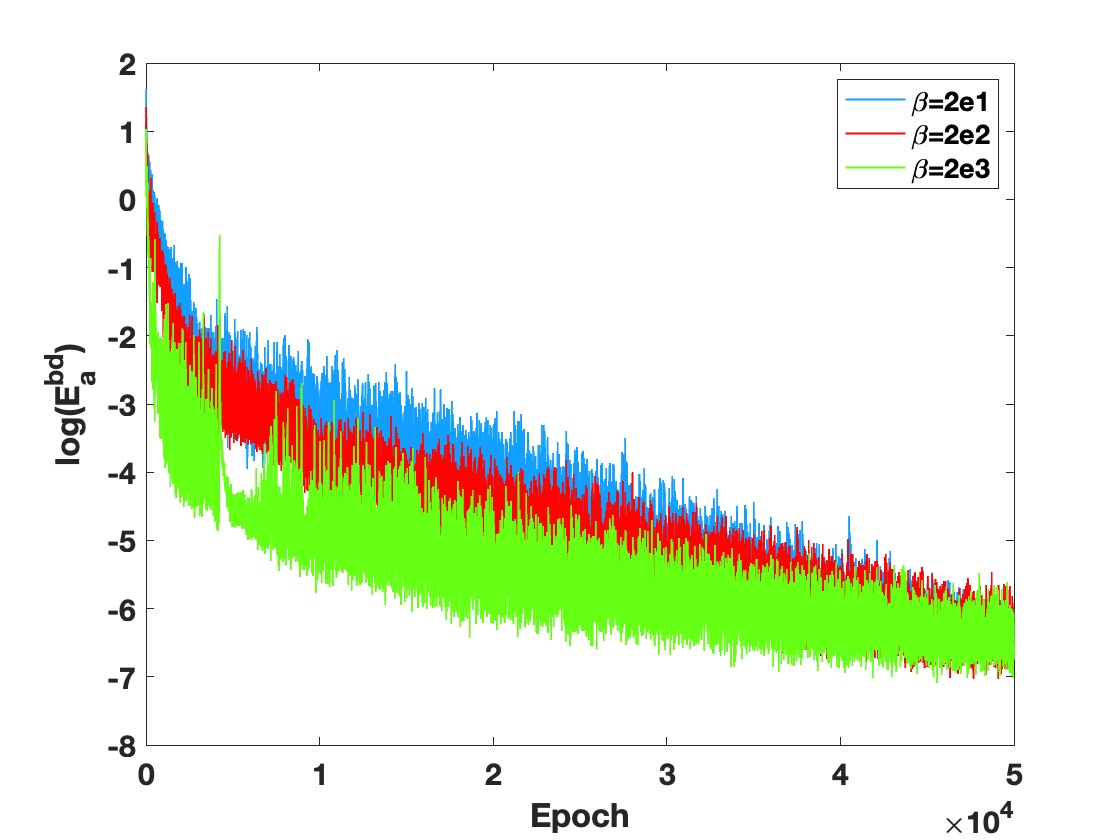}
  \end{minipage}
  \caption{The absolute error in the domain (left) and on the boundary (right) of ALDL.}
  \label{fig:2d_lineig_al}
\end{figure}

We present the numerical results by the ALDL method with different parameters $\beta$ in Figure \ref{fig:2d_lineig_al} and Table \ref{tab:2d_lineig_al}. We can see that the numerical error is approximately  $10^{-3}$  for the eigenfunctions and $10^{-4}$ for eigenvalues. Moreover, the approximation accuracy and convergence rate are robust with respect to different $\beta.$

\begin{table}[H]
\small
\begin{minipage}[H]{0.5\linewidth}
  \centering
  \begin{tabular}{cccc}
    \hline
    $\beta$ &$\mathcal{E}_{a}^{in}$ &$\mathcal{E}_{a}^{bd}$ &$|\rho_{dl}-\rho|/\rho$ \\
    \hline
     2e+1   & 5.1201e-3 &  1.6031e-3 & 2.0678e-4 \\
     2e+2   & 7.7265e-3 &  1.4665e-3 & 7.3460e-4 \\
     2e+3   & 6.9291e-3 & 1.4485e-3 &  2.1878e-4 \\
    \hline
  \end{tabular}
  \caption{ The final absolute errors of ALDL. }
  \label{tab:2d_lineig_al}
\end{minipage}
\begin{minipage}[H]{0.5\linewidth}
  \centering
  \begin{tabular}{cccc}
  \hline
    $\beta$ &$\mathcal{E}_{a}^{in}$ &$\mathcal{E}_{a}^{bd}$ &$|\rho_{dl}-\rho|/\rho$ \\
    \hline
     2e+2  & 5.7031e-2  &  6.0427e-2 & 7.4183e-2 \\
     2e+3  &  9.9260e-3 &  7.0166e-3 & 7.0739e-3 \\
     2e+4  &  1.8037e-2 &  1.0109e-3 & 9.7971e-4 \\
    \hline
  \end{tabular}
  \caption{ The final absolute errors of PMDL. }
  \label{tab:2d_lineig_pm}
\end{minipage}
\end{table}

\begin{figure}[H]
  \begin{minipage}[H]{0.5\linewidth}
    \centering
    \includegraphics[scale = 0.14]{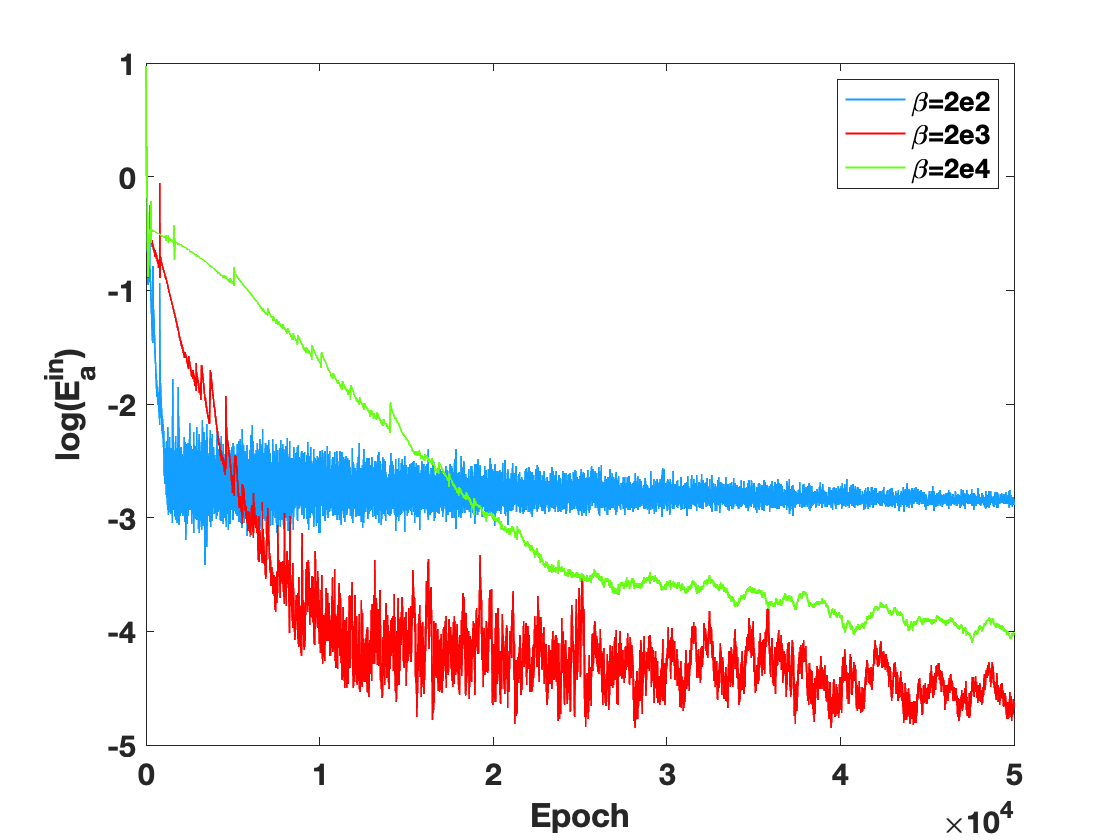}
  \end{minipage}
  \begin{minipage}[H]{0.5\linewidth}
    \centering
    \includegraphics[scale = 0.14]{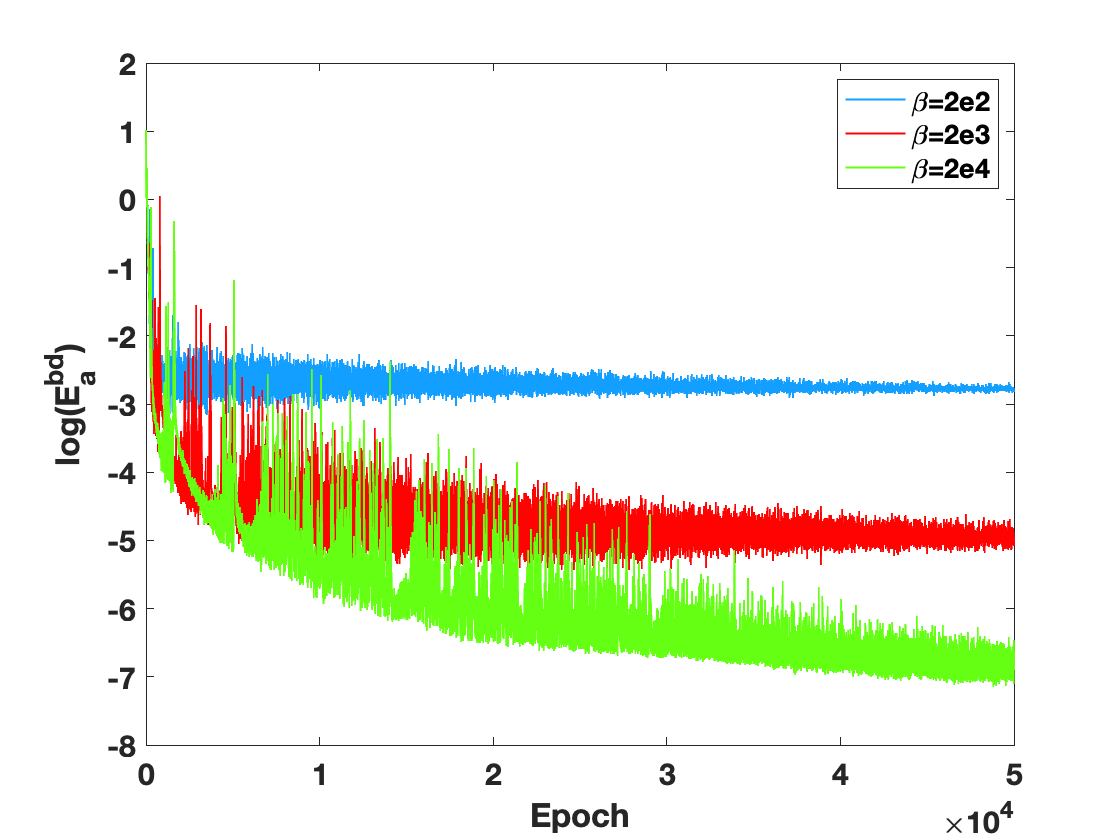}
  \end{minipage}
  \caption{ The absolute error in the domain (left) and on the boundary (right) of PMDL.}
  \label{fig:2d_lineig_pm}
\end{figure}

We also show the associated results by the PMDL method in Figure \ref{fig:2d_lineig_pm} and Table \ref{tab:2d_lineig_pm}, from which we may conclude that the accuracy and convergent speed of the PMDL method are relatively sensitive to the parameter $\beta.$ 


{We conclude from Table \ref{tab:2d_lineig_al} and Table \ref{tab:2d_lineig_pm} that the ALDL method has 2-20 times advantages over the PMDL method in terms of the accuracy of eigenfunctions or eigenvalues. Similar to the PDEs, we compare the numerical performance of the ALDL method and the SGDA method in Table \ref{tab:2d_lin_com}, from which we find that the approximations of the former are twice as accurate as that of the latter and the ALDL takes $30\%$ less time. }

\begin{table}[H]
\small
  \centering
  \begin{tabular}{cccccc}
  \hline
    Method &$\beta$ &$\mathcal{E}_{a}^{in}$ &$\mathcal{E}_{a}^{bd}$ &$|\rho_{dl}-\rho|/\rho$ &time(s) \\
    \hline
     ALDL & 2e+2  &  7.7265e-3 &  1.4665e-3 &  7.3460e-4 &  1214.60\\
     SGDA & 2e+2  &  7.6285e-3 &  2.1653e-3 &  1.0547e-3 &  1695.77\\
     ALDL & 2e+3  &  6.9291e-3 &  1.4485e-3 &  2.1878e-4 &  1008.00\\
     SGDA & 2e+3  &  9.6389e-3 &  3.0370e-3 &  1.1833e-3 &  1652.99\\
    \hline
  \end{tabular}
  \caption{ The comparation of different methods. }
  \label{tab:2d_lin_com}
\end{table}


\subsubsection{The 3D case}

\begin{figure}[H]
  \begin{minipage}[H]{0.5\linewidth}
    \centering
    \includegraphics[scale = 0.14]{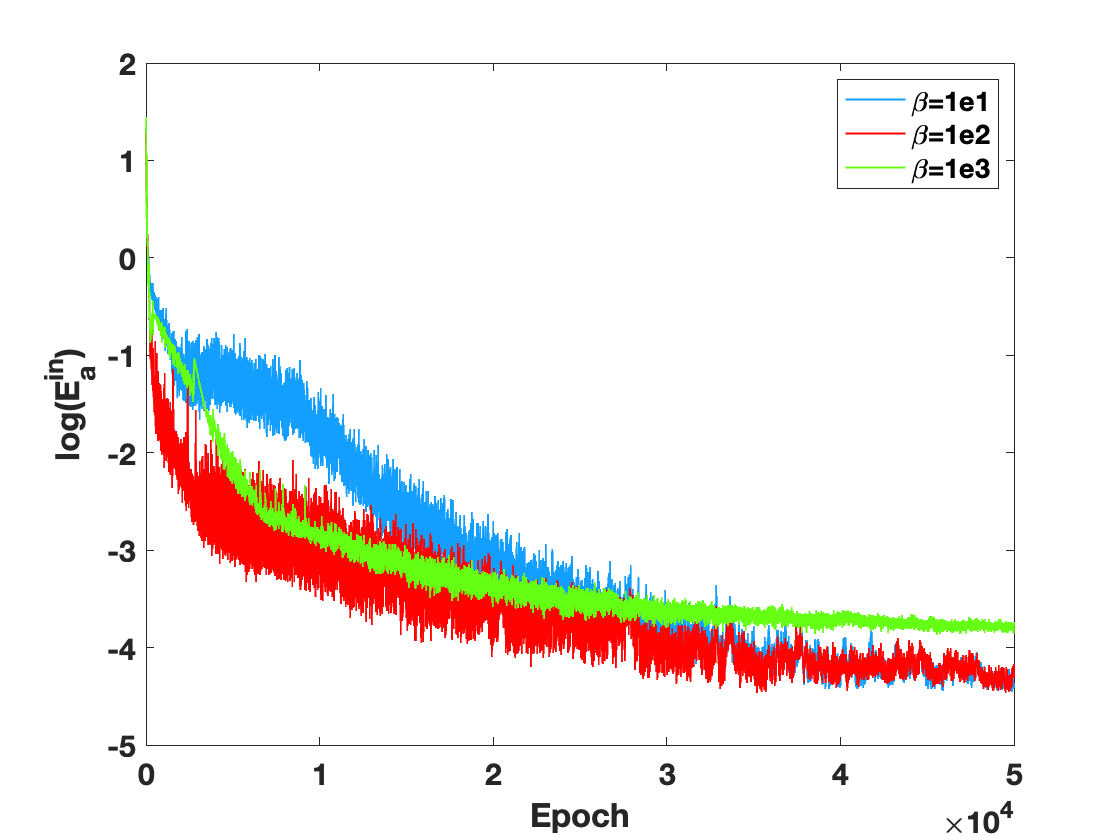}
  \end{minipage}
  \begin{minipage}[H]{0.5\linewidth}
    \centering
    \includegraphics[scale = 0.14]{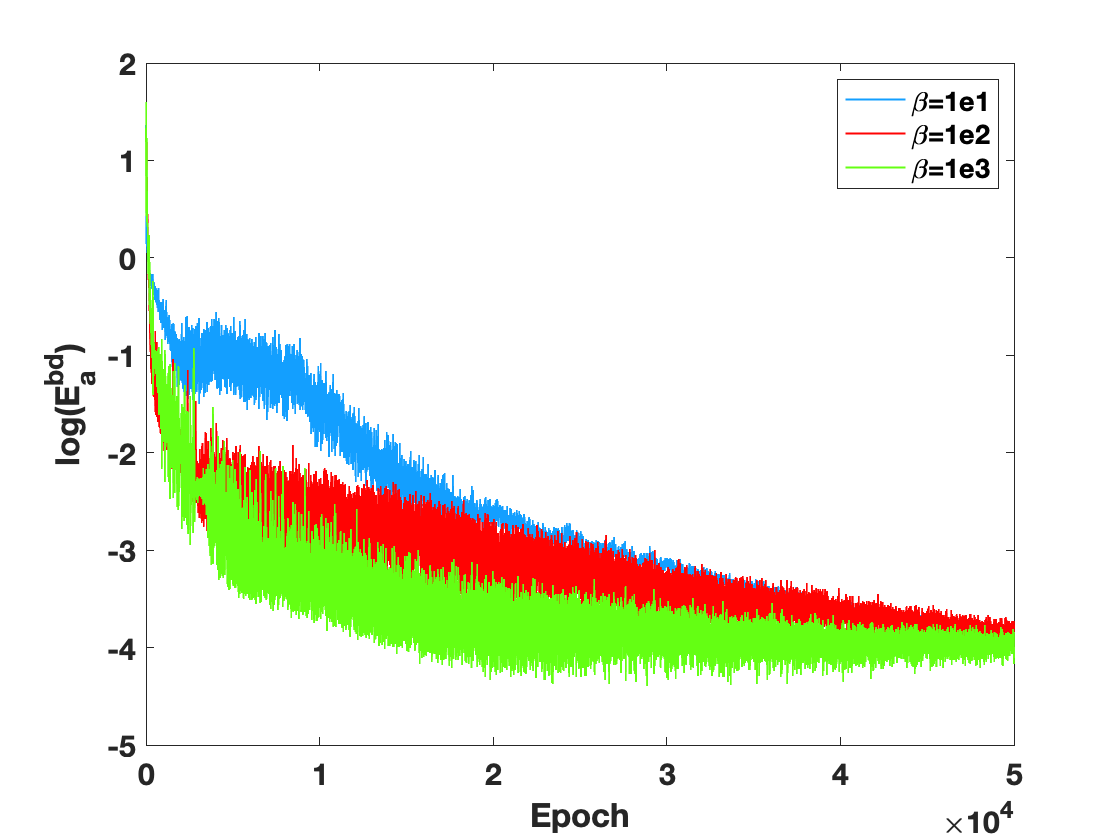}
  \end{minipage}
  \caption{ The absolute error in the domain (left) and on the boundary (right) of ALDL.}
  \label{fig:3d_lineig_al}
\end{figure}

\begin{figure}[H]
  \begin{minipage}[H]{0.5\linewidth}
    \centering
    \includegraphics[scale = 0.14]{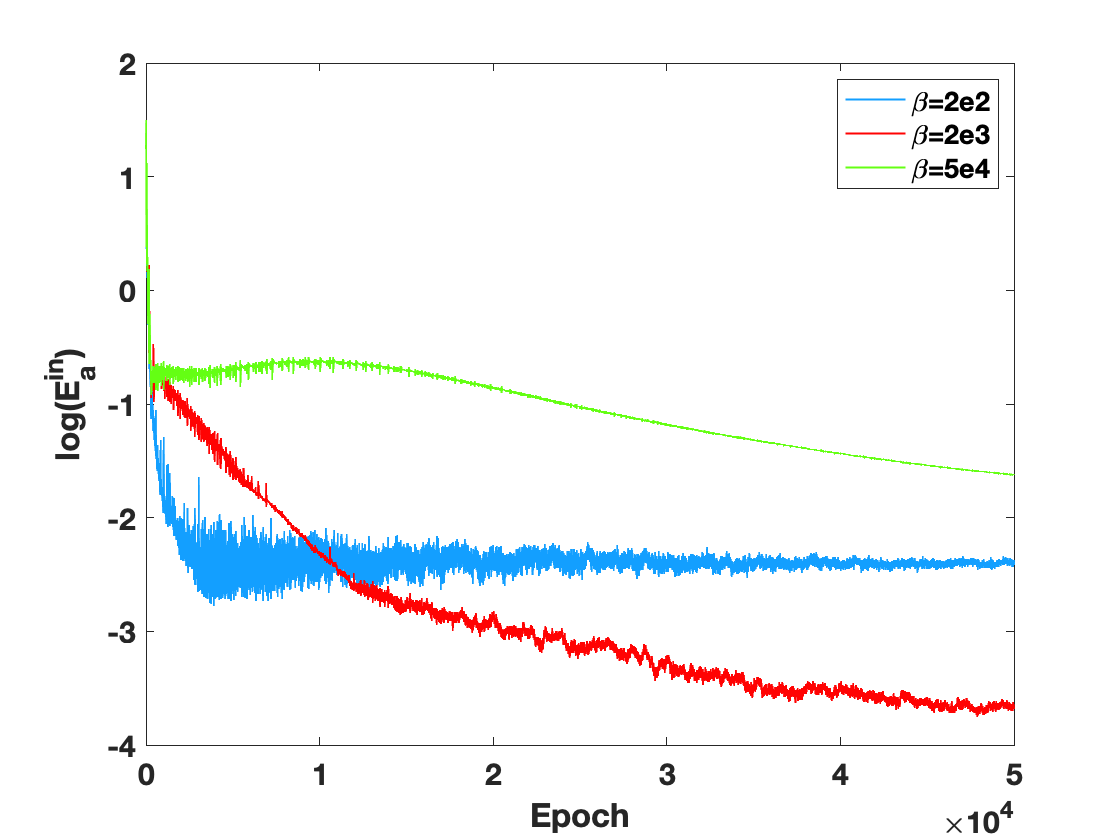}
  \end{minipage}
  \begin{minipage}[H]{0.5\linewidth}
    \centering
    \includegraphics[scale = 0.14]{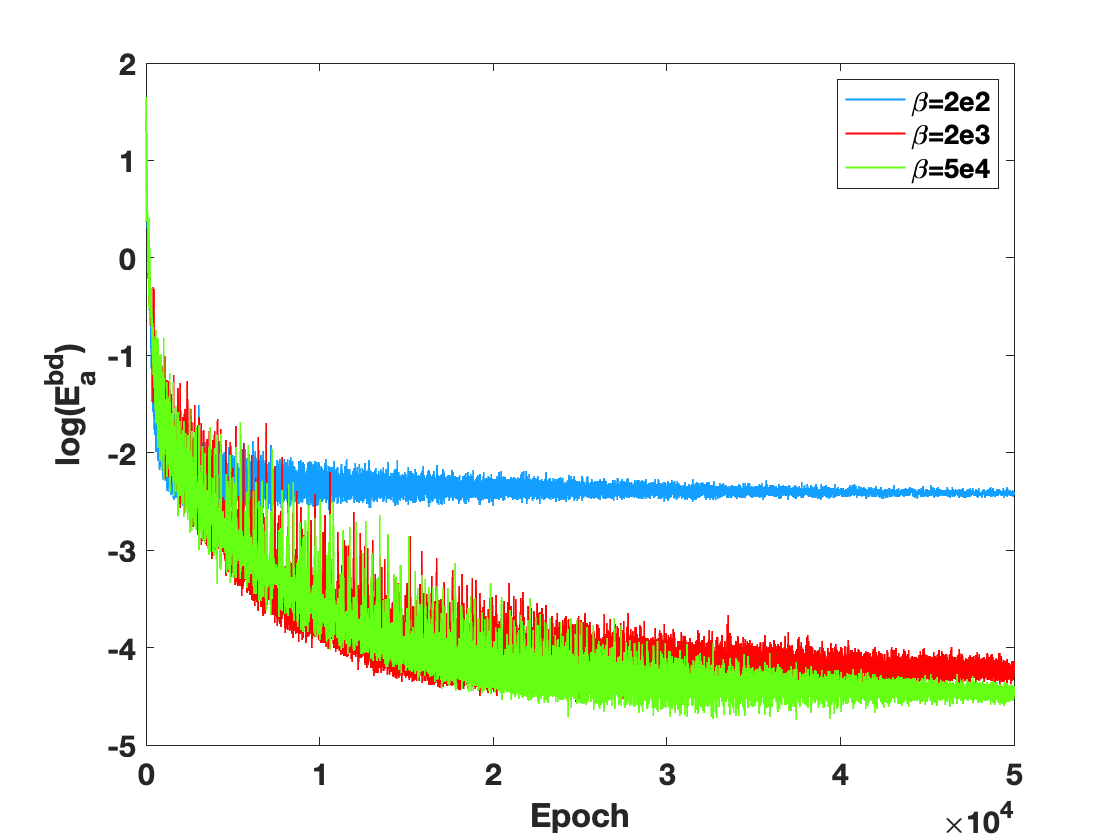}
  \end{minipage}
  \caption{ The absolute error in the domain (left) and on the boundary (right) of PMDL.}
  \label{fig:3d_lineig_pm}
\end{figure}

We have also performed the experiments for the three dimensional case, and the results are shown in Figure \ref{fig:3d_lineig_al}, Figure \ref{fig:3d_lineig_pm}, and Table \ref{tab:3d_lineig_al} - Table \ref{tab:3d_lineig_pm}. Similar to the two dimensional case, the ALDL method admits a better numerical performance in view of {the} accuracy. 
{Especially, we noticed that when the parameter $\beta$ of the PMDL method is small, the eigenvalues obtained by the ALDL method at least are more $100$ times accurate than those obtained by the PMDL method. When the parameter $\beta$ of the PMDL method is large, the eigenfunctions obtained by the ALDL method at least are more $10$ times accurate than those obtained by the PMDL method.}

\begin{table}[H]
\small
\begin{minipage}[H]{0.5\linewidth}
  \centering
  \begin{tabular}{cccc}
    \hline
    $\beta$ &$\mathcal{E}_{a}^{in}$ &$\mathcal{E}_{a}^{bd}$ &$|\rho_{dl}-\rho|/\rho$ \\
    \hline
     1e+1   &  1.2117e-2 &  2.1592e-2 &  1.4851e-3 \\
     1e+2   &  1.3606e-2 &  2.3021e-2 &  4.2040e-4 \\
     1e+3   &  2.3689e-2 &  1.7151e-2 &  1.7392e-3 \\
    \hline
  \end{tabular}
  \caption{ The final absolute errors of ALDL. }
  \label{tab:3d_lineig_al}
\end{minipage}
\begin{minipage}[H]{0.5\linewidth}
  \centering
  \begin{tabular}{cccc}
  \hline
    $\beta$ &$\mathcal{E}_{a}^{in}$ &$\mathcal{E}_{a}^{bd}$ &$|\rho_{dl}-\rho|/\rho$ \\
    \hline
     2e+2  &  9.1493e-2 &  8.9782e-2 &  7.6310e-2 \\
     2e+3  &  2.5106e-2 &  1.5402e-2 &  7.2753e-3 \\
     5e+4  &  1.9755e-1 &  1.1124e-2 &  9.5417e-3 \\
    \hline
  \end{tabular}
  \caption{ The final absolute errors of PMDL. }
  \label{tab:3d_lineig_pm}
\end{minipage}
\end{table}

\subsection{A nonlinear eigenvalue problem}

Finally, we consider the nonlinear Schr{\"o}dinger eigenvalue problem
\[
  \begin{cases}
    - \Delta u + Vu + u^3 = \rho u \quad   &\mbox{\rm in}  \ \Omega,
    \\
    \quad u = 0 \quad &\mbox{\rm on} \ \Gamma,
    \\
    \|u\|_{0,\Omega} = 1,
  \end{cases}
\]
where $\Omega = (0,1)^d$ for $d=2$ and $3$, and $V = \sum_{i=1}^d x_i^2$. Since the true solutions $\rho$ and $u$ are unavailable, we shall {numerically} compute the reference solutions $\rho_{ref}$ and $u_{ref}.$ More precisely, we construct a DNN reference solution $u_{ref}$ by:
\[
  u_{ref}(\boldsymbol{x};\boldsymbol{\theta}) = \ell(\boldsymbol{x})\psi(\boldsymbol{x};\boldsymbol{\theta}),
\]
where $\psi(\boldsymbol{x};\boldsymbol{\theta})$ is a ResNet function with width 50 and depth 6, and  $$\ell(\boldsymbol{x}) = \prod_{i=1}^d x_i (1-x_i).$$
In this way the reference solution admits an exact match on the boundary. The reference solution  $u_{ref}(\boldsymbol{x};\boldsymbol{\theta})$ is trained by the Adam optimizer with a learning rate $\eta=5e-4,$ and we set $Epoch = 100000$ with batch size $2048$ in the domain { for $2$d problem and $8192$ for $3$d problem.}

\subsubsection{The 2D case}

\begin{figure}[H]
  \begin{minipage}[H]{0.5\linewidth}{}
    \centering
    \includegraphics[scale = 0.14]{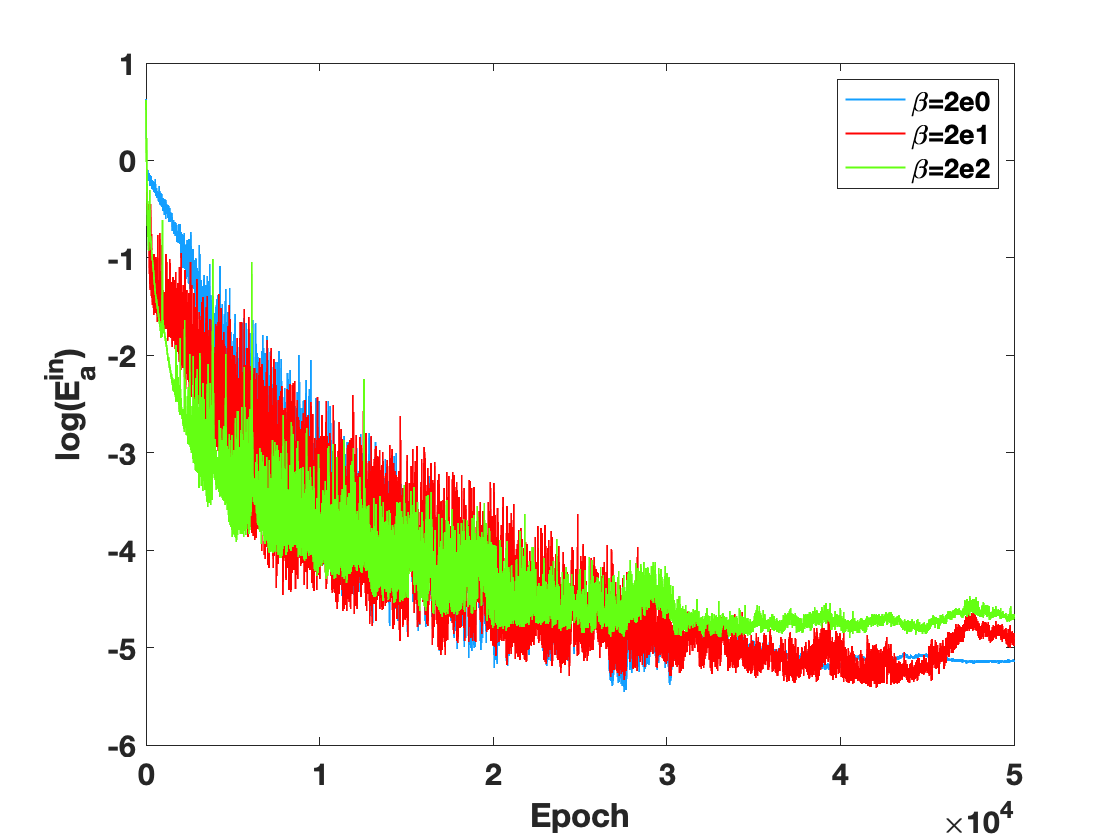}
  \end{minipage}
  \begin{minipage}[H]{0.5\linewidth}
    \centering
    \includegraphics[scale = 0.14]{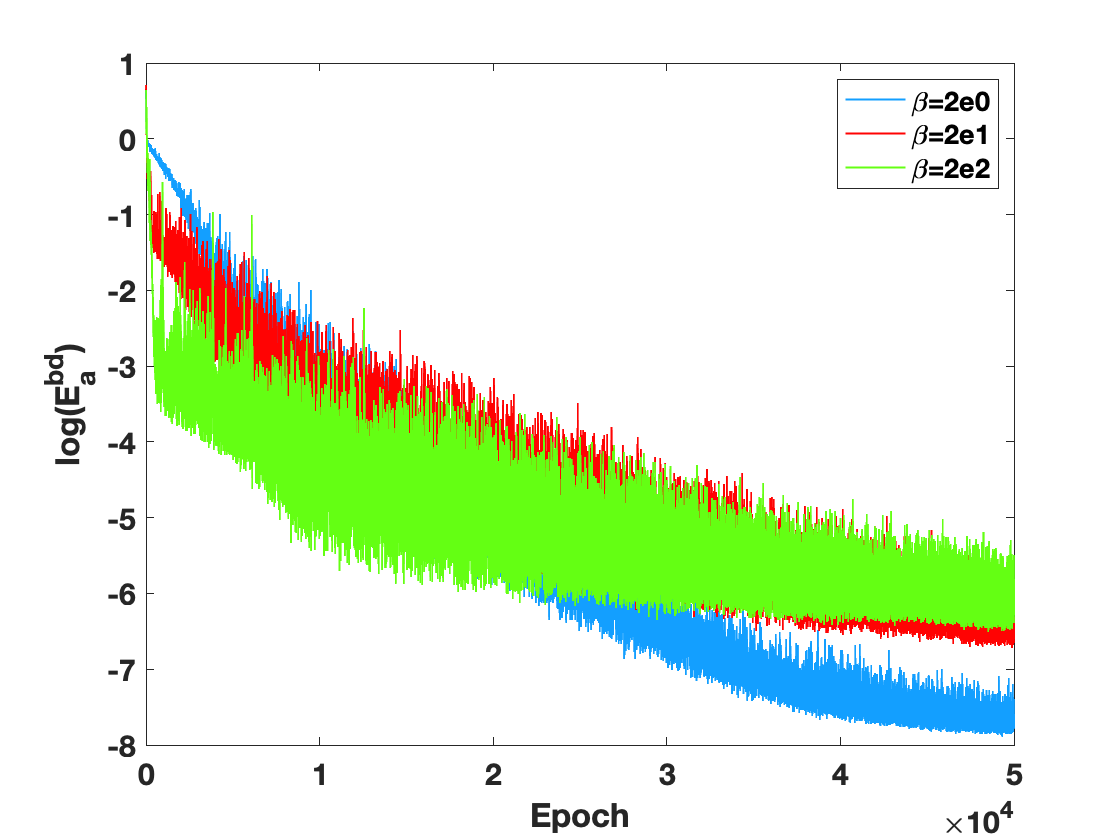}
  \end{minipage}
  \caption{ The absolute error in the domain (left) and on the boundary (right) of ALDL.}
  \label{fig:2d_noneig_al}
\end{figure}

The numerical results for the ALDL method with different parameters $\beta$ are presented in Figure \ref{fig:2d_noneig_al} and Table \ref{tab:2d_noneig_al}. Again, it is noticed that the method is insensitive to the choice of $\beta.$

\begin{table}[H]
\small
\begin{minipage}[H]{0.5\linewidth}
  \centering
  \begin{tabular}{cccc}
    \hline
    $\beta$ &$\mathcal{E}_{a}^{in}$ &$\mathcal{E}_{a}^{bd}$ &$|\rho_{dl}-\rho_{ref}|/\rho_{ref}$ \\
    \hline
     2e+0   &  5.8958e-3 &  4.8274e-4 &  3.5429e-3 \\
     2e+1   &  7.3927e-3 &  1.8059e-3 &  1.6044e-4\\
     2e+2   &  9.2868e-3 &  2.9953e-3 &  5.2254e-4  \\
    \hline
  \end{tabular}
  \caption{ The final absolute errors of ALDL. }
  \label{tab:2d_noneig_al}
\end{minipage}
\begin{minipage}[H]{0.5\linewidth}
  \centering
  \begin{tabular}{cccc}
  \hline
    $\beta$ &$\mathcal{E}_{a}^{in}$ &$\mathcal{E}_{a}^{bd}$ &$|\rho_{dl}-\rho_{ref}|/\rho_{ref}$ \\
    \hline
     2e+2  &  3.3200e-2 &  3.3776e-2 &  4.0211e-2 \\
     2e+3  &  1.9137e-2 &  4.1947e-3 &  3.5401e-3 \\
     2e+4  &  1.2269e-1 &  2.6753e-3 &  4.1881e-3 \\
    \hline
  \end{tabular}
  \caption{ The final absolute errors of PMDL. }
  \label{tab:2d_noneig_pm}
\end{minipage}
\end{table}

\begin{figure}[H]
  \begin{minipage}[H]{0.5\linewidth}
    \centering
    \includegraphics[scale = 0.14]{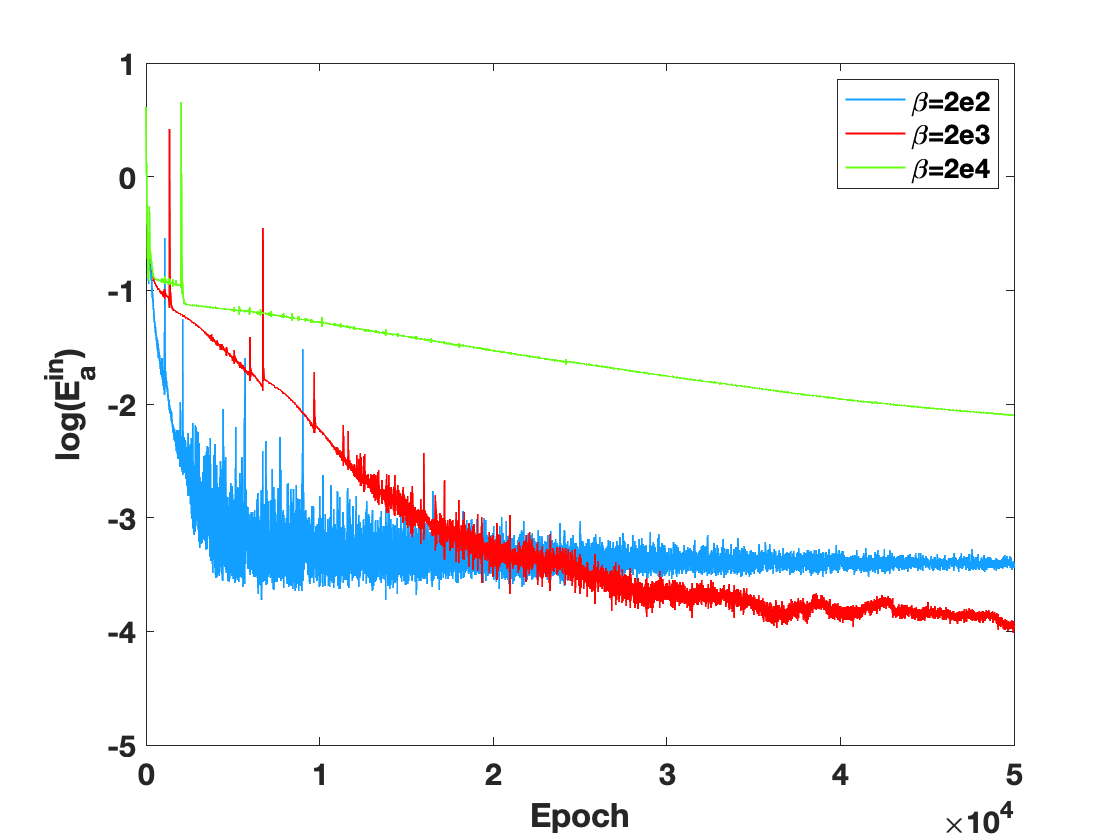}
  \end{minipage}
  \begin{minipage}[H]{0.5\linewidth}
    \centering
    \includegraphics[scale = 0.14]{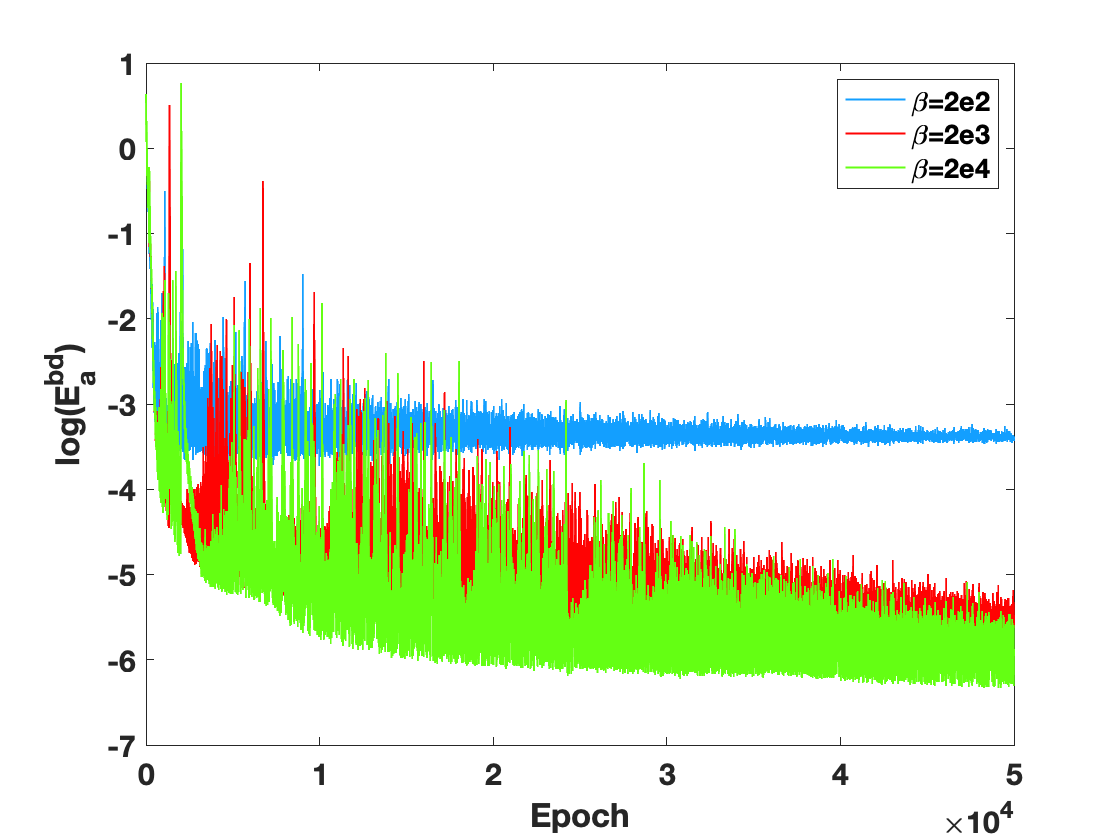}
  \end{minipage}
  \caption{ The absolute error in the domain (left) and on the boundary (right) of PMDL.}
  \label{fig:2d_noneig_pm}
\end{figure}


The associated numerical results by the PMDL method are presented in Figure \ref{fig:2d_noneig_pm} and Table \ref{tab:2d_noneig_pm}, from which we can see that the larger the $\beta$ is, the slower the converge rate seems to be. {In addition, the ALDL method has $2-20$ times advantages over the PMDL method in terms of the accuracy of eigenfunctions. As for eigenvalues, the ALDL method has more than $100$ times the accuracy advantage over the PMDL method in some cases.}


\subsubsection{The 3D case}

\begin{figure}[H]
  \begin{minipage}[H]{0.5\linewidth}
    \centering
    \includegraphics[scale = 0.14]{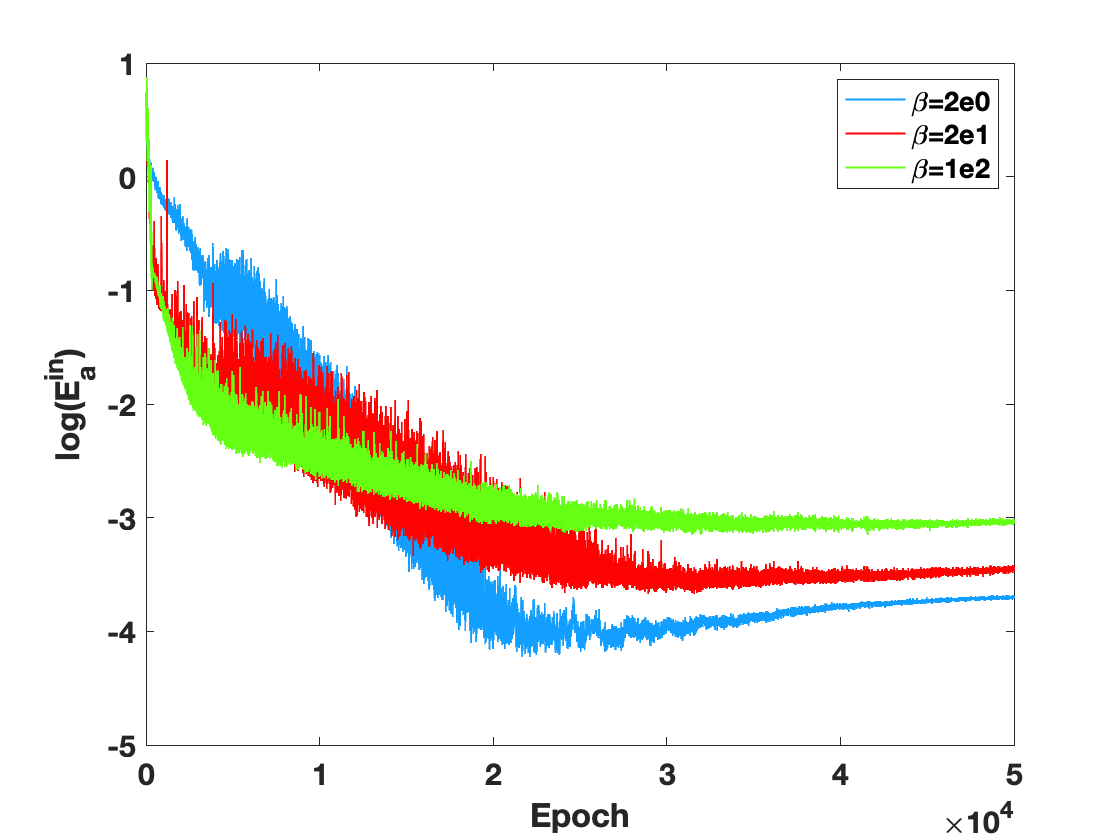}
  \end{minipage}
  \begin{minipage}[H]{0.5\linewidth}
    \centering
    \includegraphics[scale = 0.14]{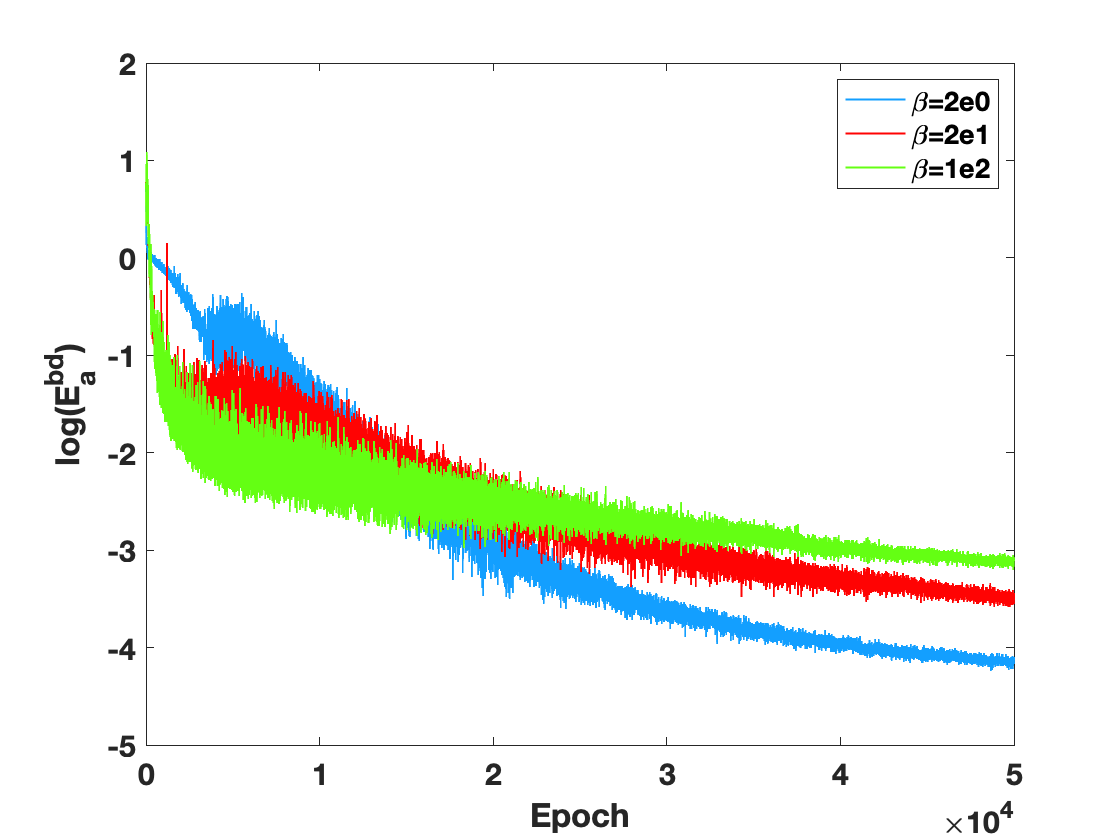}
  \end{minipage}
  \caption{ The absolute error in the domain (left) and on the boundary (right) of ALDL.}
  \label{fig:3d_noneig_al}
\end{figure}

\begin{figure}[H]
  \begin{minipage}[H]{0.5\linewidth}
    \centering
    \includegraphics[scale = 0.14]{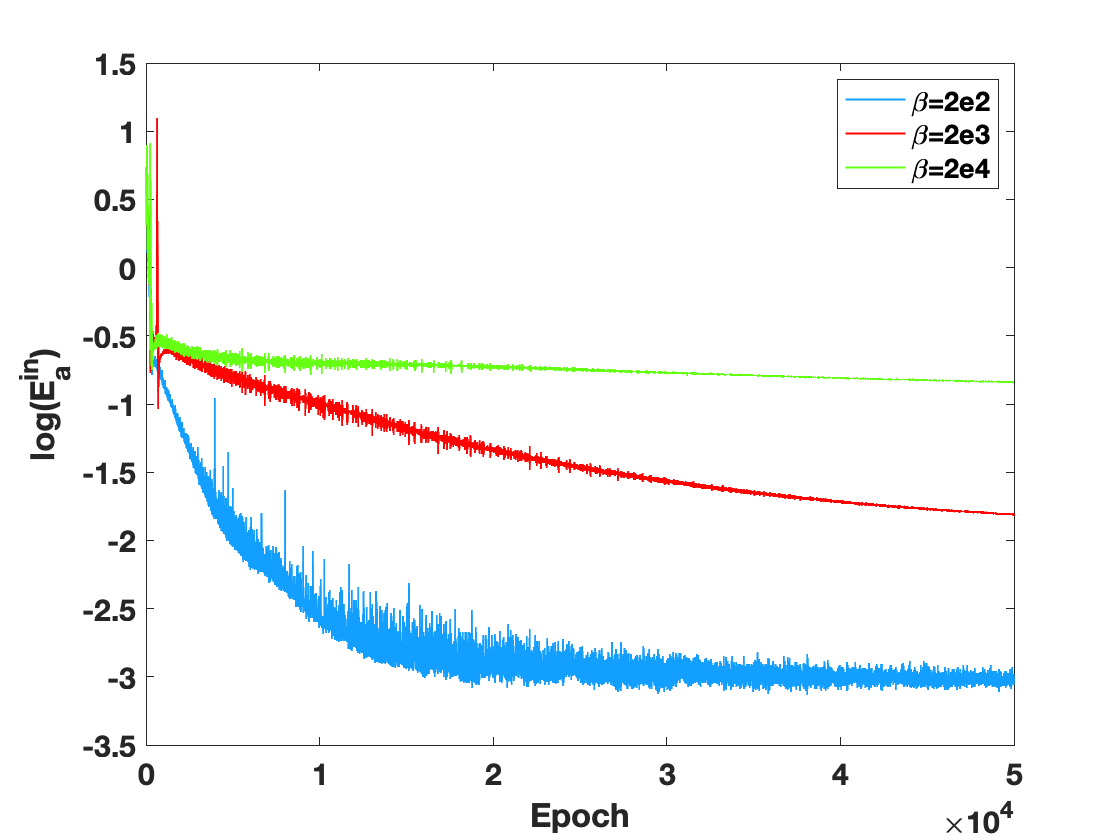}
  \end{minipage}
  \begin{minipage}[H]{0.5\linewidth}
    \centering
    \includegraphics[scale = 0.14]{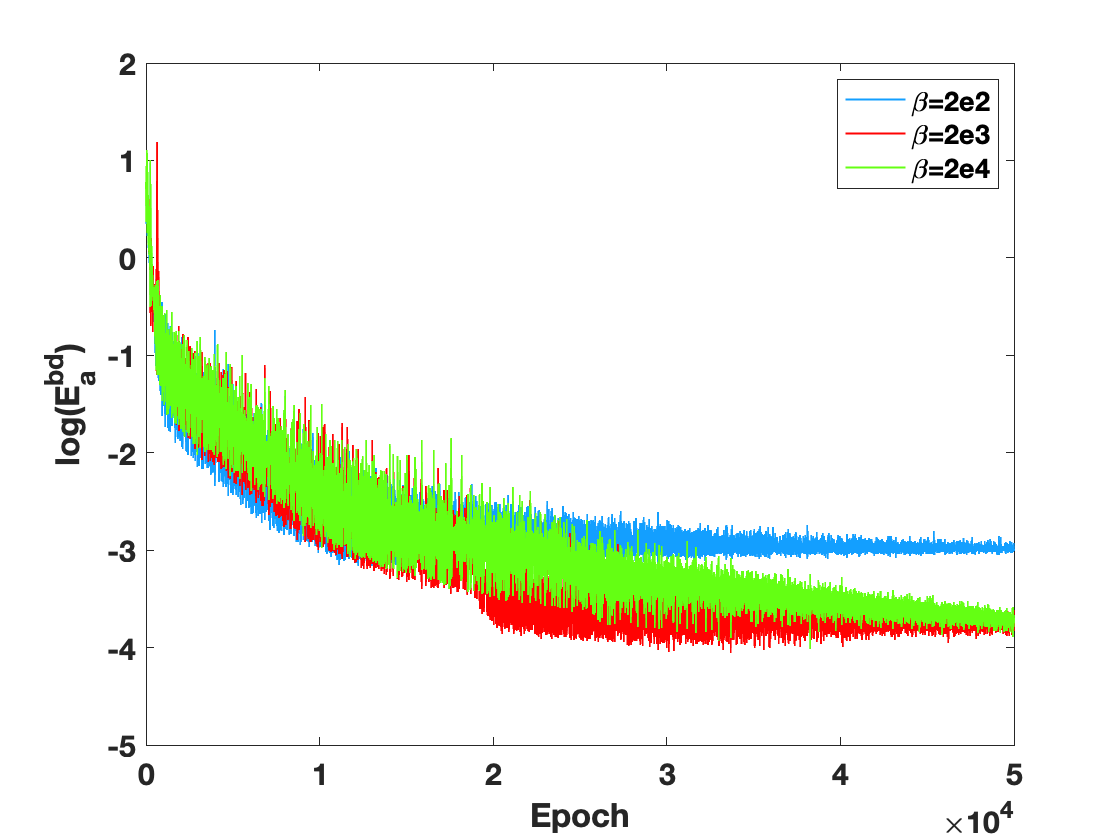}
  \end{minipage}
  \caption{ The absolute error in the domain (left) and on the boundary (right) of PMDL.}
  \label{fig:3d_noneig_pm}
\end{figure}

We have also presented the 3D simulations, and the results are shown in Figure \ref{fig:3d_noneig_al} - Figure \ref{fig:3d_noneig_pm} and Table \ref{tab:3d_noneig_al} -  Table \ref{tab:3d_noneig_pm}, and one can draw similar conclusions as in the above examples. 

\begin{table}[H]
\small
\begin{minipage}[H]{0.5\linewidth}
  \centering
  \begin{tabular}{cccc}
    \hline
    $\beta$ &$\mathcal{E}_{a}^{in}$ &$\mathcal{E}_{a}^{bd}$ &$|\rho_{dl}-\rho_{ref}|/\rho_{ref}$ \\
    \hline
     2e+0   &  2.4886e-2 &  1.6198e-2 &  3.0750e-3 \\
     2e+1   &  3.1448e-2 &  3.1834e-2 &  4.0040e-4 \\
     1e+2   &  4.7370e-2 &  4.5080e-2 &  7.2167e-4 \\
    \hline
  \end{tabular}
  \caption{The final absolute errors of ALDL. }
  \label{tab:3d_noneig_al}
\end{minipage}
\begin{minipage}[H]{0.5\linewidth}
  \centering
  \begin{tabular}{cccc}
  \hline
    $\beta$ &$\mathcal{E}_{a}^{in}$ &$\mathcal{E}_{a}^{bd}$ &$|\rho_{dl}-\rho_{ref}|/\rho_{ref}$ \\
    \hline
     2e+2  &  4.9407e-2 &  5.0387e-2 &  4.0671e-2 \\
     2e+3  &  1.6322e-1 &  2.4755e-2 &  7.3079e-3 \\
     2e+4  &  4.3183e-1 &  2.3657e-2 &  1.0364e-2 \\
    \hline
  \end{tabular}
  \caption{ The final absolute errors of PMDL. }
  \label{tab:3d_noneig_pm}
\end{minipage}
\end{table}

\section{Concluding remarks}\label{sec: con}

We have proposed an augmented Lagrangian deep learning  method for variational problems with essential boundary conditions. The approach relies on first rewriting the original problem into an equivalent minimax problem, and then expressing the primal and dual variables with two individual DNN functions. Then, the network parameters of the primal and dual variables are trained using the stochastic optimization method together with a projection technique. Applications to elliptic problems and eigenvalue problems show that the ALDL method admits many advantages over the penalty method. In our future studies, we shall extend our ALDL approach to time dependent problems with complex solution structures.

\section*{Acknowledgments}
J. Huang is partially supported by the National Key Research and Development Project (Grant No. 2020YFA0709800), NSFC (Grant No. 12071289) and Shanghai Municipal Science and Technology Major Project (2021SHZDZX0102). 
T. Zhou is supported by the National Key R\&D Program of China (2020YFA0712000), NSFC (under grant numbers 11822111, 11688101), the science challenge project (No. TZ2018001), and youth innovation promotion association (CAS). 

\bibliographystyle{plain}

\end{document}